\documentclass[11pt,runningheads,amsmath,amssymb, eucal]{ZeClasse}

\usepackage{latexsym}
\usepackage{fullpage}
\usepackage{amsfonts}
\usepackage{epsfig}
\input xy
\xyoption{all}

\makeatletter 
\newcommand\mysection{\@startsection {section}{1}{\z@}%
                                    {0ex \@plus -1ex \@minus -.2ex}%
                                    {\thepec}%
                                    {\normalfont\normalsize\bfseries}}
\newcommand\mysubsection{\@startsection{subsection}{2}{\z@}%
                                      {0ex\@plus -1ex \@minus -.2ex}%
                                      {-1em}%
                                      {\normalfont\normalsize\it}}
\newcommand\mysubsubsection{\@startsection{subsubsection}{3}{\z@}%
                                      {0ex\@plus -1ex \@minus -.2ex}%
                                      {-1em}%
                                      {\normalfont\normalsize\it}}
\makeatother 

\newtheorem{thm}{Theorem}
\newtheorem{prop}[theorem]{Proposition}

\newcommand{\prv}{\ensuremath{\mathbf{P}V}}
\newcommand{\carre}{\hspace{\fill} $\Box$}

\begin{document}

\input amssym.def
\input amssym.tex

\def\page#1{\leaders\hbox to 5mm{.}\hfill \rlap{\hbox to 2mm{\hfill #1}}\par}

\font\tengoth=eufm10 at 12pt
\font\sevengoth=eufm7
\font\fivegoth=eufm5
\newfam\gothfam
\textfont\gothfam=\tengoth
\scriptfont\gothfam=\sevengoth
\scriptscriptfont\gothfam=\fivegoth
\def\goth{\fam\gothfam\tengoth}

\def\thepec{2mmplus1mmminus1mm} 
\def\thegec{4mmplus1mmminus1mm} 
\def\pec{\vskip\thepec}
\def\gec{\vskip\thegec}
\def\pech{\hskip3mmplus1mmminus1mm}
\def\gech{\hskip6mmplus1mmminus1mm}
\def\lc{\hbox{\rm [loc. cit.]}~}

\def\og{\leavevmode\raise.3ex\hbox{$\scriptscriptstyle\langle\!\langle\,$}}
\def\fg{\leavevmode\raise.3ex\hbox{$\scriptscriptstyle\,\rangle\!\rangle\ $}}

\def\limind{\mathop{\oalign{lim\cr \hidewidth$\longrightarrow$\hidewidth}}}
\def\limproj{\mathop{\oalign{lim\cr \hidewidth$\longleftarrow$\hidewidth}}}
\def\tch#1{#1\mkern2.5mu\check{}}



\title*{Topological simplicity, commensurator super-rigidity and 
non-linearities of Kac-Moody
groups}
\titlerunning{Topological simplicity and non-linearities of Kac-Moody groups}

\author{{\bf Appendix by P. Bonvin: Strong boundaries and 
commensurator super-rigidity}
\\
\vskip 6mm Bertrand R\'EMY}
\authorrunning{Bertrand R\'EMY}

\institute{
Institut Fourier -- UMR CNRS 5582\\
Universit\'e de Grenoble 1 -- Joseph Fourier\\
100, rue des maths -- BP 74\\
38402 Saint-Martin-d'H\`eres Cedex -- France\\
{\tt bremy@fourier.ujf-grenoble.fr}}

\maketitle

\gec

\centerline{Pr\'epublication de l'Institut Fourier n$^0$ {\bf 590} (2003)}
\centerline{http:$/\!\!/$www-fourier.ujf-grenoble.fr/prepublications.html}

\gec

{\footnotesize
{\sc Abstract}.---~
We provide new arguments to see topological Kac-Moody groups as 
generalized semisimple groups
over local fields: they are products of topologically simple groups 
and their Iwahori
subgroups are the normalizers of the pro-$p$ Sylow subgroups.
We use a dynamical characterization of parabolic subgroups to prove 
that some countable
Kac-Moody groups with Fuchsian buildings are  not linear.
We show for this that the linearity of a countable Kac-Moody group 
implies the existence of a
closed  embedding of the corresponding topological group in a 
non-Archimedean simple Lie
group, thanks to a commensurator super-rigidity theorem proved in the 
Appendix by P. Bonvin.}

\pec

{\footnotesize

}

\vfill

\centerline{\bf Contents}

\pec

{\bf Introduction}
\page{2}

\pec

{\bf 1. Structure theorem}

1.A Topological Kac-Moody groups
\page{4}

1.B Refined Tits system and virtual pro-$p$-ness of parahoric subgroups
\page{6}

1.C Lattices and generalized arithmeticity
\page{7}

\pec

{\bf 2. Topological simplicity theorem}

2.A Topological simplicity
\page{9}

2.B Discussion of the proof for countable groups
\page{11}

2.C Closures of Levi factors
\page{12}

\pec

{\bf 3. Embedding theorem}

3.A Semisimple Zariski closure and injectivity
\page{15}

3.B Unbounded image and continuous extension
\page{16}

3.C Embedding of vertices and closed image
\page{16}

\pec

{\bf 4. Some concrete non-linear examples}

4.A Groups with right-angled Fuchsian buildings
\page{17}

4.B Dynamics and parabolics
\page{21}

4.C Non-linearity in equal characteristic
\page{23}

\pec

{\bf Appendix, by P. Bonvin: Strong boundaries and commensurator 
super-rigidity}
\page{27}

\eject

{\bf Introduction}

\pec

This paper contains two kinds of results, according to which 
definition of Kac-Moody groups is
adopted.
The main goal is to prove non-linearities of countable Kac-Moody groups as
defined by generators and relations by J. Tits \cite{Tit87}.
The strategy we follow, based on continuous 
extensions of abstract group
homomorphisms, leads us to prove structure results on topological 
Kac-Moody groups as previously
defined in \cite{RR02}.
We first give the statements; we will provide below details and 
motivations for each
theorem individually, and then explain the relationship between them.
The first theorem basically says that infinitely many countable 
Kac-Moody groups
are not linear over any field (Theorem \ref{sss - RA not 
linear}).  Its proof uses
super-rigidity and dynamical properties of generalized parabolic subgroups.

\pec

{\sc Theorem}.---~\it
Let $\Lambda$ be a generic countable Kac-Moody group over a finite 
field with right-angled
Fuchsian buildings.
Then $\Lambda$ is not linear over any field.
Moreover for each prime number $p$, there are infinitely many such 
non-linear groups defined over
a finite field of characteristic $p$.
\rm\pec

The second theorem says that topological Kac-Moody groups should  
be seen as generalizations
of semisimple groups over local fields of positive characteristic.
We learnt from O. Mathieu that R.V. Moody has proved some topological 
simplicity results in Kac-Moody theory; unfortunately we couldn't find any 
written version of this work. 
The statement below roughly sums up Theorem \ref{sss - top simple} 
and Proposition
\ref{sss - Sylow}.

\pec

{\sc Theorem}.---~\it
{\rm (i)~} A topological Kac-Moody group over a finite field is the 
direct product of
topologically simple groups, with one factor for each connected 
component of its Dynkin
diagram.

{\rm (ii)~} The Iwahori subgroups, i.e. the chamber fixators for the 
natural action on the
building of the group, are characterized as the normalizers of the 
pro-$p$ Sylow subgroups.
\rm\pec

The motivation for both results is based on an analogy with 
well-known classical cases. 
Namely, Kac-Moody theory is usually introduced via the specific affine type.
For Lie algebras, this case corresponds to semisimple Lie algebras 
tensorised by
Laurent polynomials over the groundfield \cite{Kac90}.
For countable Kac-Moody groups over finite fields, it corresponds to
$\{ 0;\infty \}$-arithmetic subgroups of semisimple groups ${\bf G}$ 
over function fields.
That general countable Kac-Moody groups over finite fields provide a 
generalization of
arithmetic groups was observed some years ago.
The main argument is that there is a natural diagonal action of such 
a group $\Lambda$ on the
product of two isomorphic buildings, and when the groundfield is 
large enough, $\Lambda$ is a
lattice of the product of (the automorphism groups of) the buildings 
\cite{CG99},
\cite{Rem98}.
In the affine case, when $\Lambda={\bf G}({\bf F}_q[t,t^{-1}])$, the 
two buildings are the
Bruhat-Tits buildings of ${\bf G}$ over the completions ${\bf 
F}_q(\!(t)\!)$ and
${\bf F}_q(\!(t^{-1})\!)$.

\pec

In the general case, it is obvious that some buildings are new, e.g. 
because they admit
natural ${\rm CAT}(-1)$ metrics.
Still, on the group side, we could imagine a situation where a 
well-known discrete group
operates on exotic geometries.
A first step to prove that this is not the case was to prove that
for a countable Kac-Moody group over ${\bf F}_q$, the only possible 
linearity is
over a field of the same characteristic $p$ \cite{Rem01}.
The first quoted theorem above shows that there exist countable 
Kac-Moody groups over
finite fields which are not linear, even in equal characteristic.
Of course, the affine example of arithmetic groups shows that the 
answer to the linearity
question cannot be stated in the equal characteristic case as for 
inequal characteristics.
Arithmetic groups also make expect that the proof of the former case 
is not as easy as in
the latter case.

\pec

This is indeed the case, but the proof is fruitful since it gives 
structure results
for topological Kac-Moody groups.
Such a group is defined in \cite{RR02} as the closure of the 
non-discrete action of a
countable Kac-Moody group on only one of the two buildings, 
forgetting the other.
In [loc. cit.] it was proved that these groups satisfy the axioms of 
a sharp refinement of
Tits systems and that the parahoric subgroups, i.e. the facet 
fixators, are virtually pro-$p$.
Once again, these results are {\it a posteriori~} not surprising when 
considering the affine
case: then the topological Kac-Moody groups correspond to groups of the form
${\bf G}({\bf F}_q(\!(t)\!))$ with ${\bf G}$ as above.
The second quoted theorem makes deeper the analogy between 
topological Kac-Moody groups and
semisimple groups over local fields of positive characteristic -- see 
\ref{sss - refined TS}
for further arguments.

\pec

We can now explain the connection between the two results.
The general strategy of the proof consists in strengthening the 
linearity assumption for a
countable Kac-Moody group to obtain a closed embedding of the 
corresponding topological
Kac-Moody group into a simple non-Archimedean Lie group.
The main tool is a very general commensurator super-rigidity, stated 
by M. Burger \cite{Bur95}
according to ideas of G. Margulis \cite{Mar90} and proved in the 
Appendix by P. Bonvin.
We found the idea to use super-rigidity in order to disprove linearity in a 
paper written by A. Lubotzky, Sh. Mozes and R.J. Zimmer \cite{LMZ94}. 
In Sect. \ref{s - Embedding theorem}, we prove the following dichotomy. 

\pec

{\sc Theorem}.---~\it
Let $\Lambda$ be a countable Kac-Moody group with connected Dynkin 
diagram and large enough
groundfield.
Under mild assumptions, either $\Lambda$ is not linear or the 
corresponding topological
group is a closed subgroup of a simple non-Archimedean Lie group $G$, 
with equivariant
embedding of the vertices of the Kac-Moody building to the vertices 
of the Bruhat-Tits
building of $G$.
\rm\pec

To prove this result we need further structure results about closures 
of Levi factors:
these groups have a Tits sub-system, are virtually products of 
topologically simple factors and
their buildings naturally appear, via an inessential geometric realization, 
in the building of
the ambient Kac-Moody group (\ref{sss - Levi factors}).

\pec

This criterion is used to obtain complete non-linearities by proving
that some topological Kac-Moody groups are not closed subgroups of Lie groups.
We concentrate at this stage on Kac-Moody groups with right-angled 
Fuchsian buildings, the
most well-understood hyperbolic buildings.
We define generalized parabolic subgroups as boundary point 
stabilizers, in complete analogy
with the symmetric space or Bruhat-Tits building case.
The dynamical approach to these groups is the last argument.
We use G. Prasad's work showing that to each suitable element of a
non-Archimedean semisimple group is attached a proper parabolic 
subgroup \cite{Pra77}.
These results, as well as their analogues for groups with 
right-angled Fuchsian buildings,
enable us to exploit the incompatibility between the geometries of 
${\rm CAT}(-1)$ Kac-Moody
and Euclidean Bruhat-Tits buildings (\ref{sss - RA not linear}).
This is then geometrically explained in the framework of dynamics of 
translations on
group-theoretical (or Furstenberg) compactifications of buildings (\ref{sss -
compactifications}).

\pec

Discrete group-theoretists may object that a direct way to disprove 
linearity of
finitely generated groups is to disprove their residual finiteness.
We chose the opposite point of view, at least because it led us to 
prove structure results
for a new class of topological groups.
Now that we know that some Kac-Moody groups are
non-linear, we at last feel comfortable to ask ourselves whether some 
of them are
non-residually finite, or even simple.  For this question, the work 
by M. Burger and Sh. Mozes
\cite{BM00} on lattices of products of trees may be relevant.
Moreover the work by R. Pink on compact subgroups of non-Archimedean 
Lie groups \cite{Pin98}
may lead to a complete answer to the linearity question of countable 
Kac-Moody groups, purely
in terms of the Dynkin diagram (and the size of the finite groundfield).

\pec

This article is written as follows.
In Sect. 1, we provide some references on Kac-Moody groups as defined 
by J. Tits, and recall
some facts, mainly of combinatorial nature, on topological Kac-Moody 
groups, as well as the
analogy with algebraic groups.
The new result here is the intrinsic characterization of Iwahori subgroups.
The discrete group point of view is also briefly discussed.
In Sect. 2, we prove the topological simplicity theorem, based on 
Tits system and pro-$p$
group arguments.
We explain why the proof doesn't work for countable groups, but also 
why some other simplicity
results in this case are not excluded.
In Sect. 3, we use the commensurator super-rigidity theorem proved in 
the Appendix by P.
Bonvin to strengthen the linearity assumption as in the last theorem 
above: the linearity of a
countable Kac-Moody groups leads to a closed embedding of the 
corresponding topological group
in an algebraic one.
In Sect. 4, we define and analyze generalized parabolic subgroups of 
topological Kac-Moody
groups with right-angled Fuchsian buildings.
After proving a decomposition involving a generalized unipotent 
radical, we turn to the
dynamical study of these groups, in the spirit of G. Prasad's work.
We conclude with the proof of the non-linearity theorem, followed by 
a geometric explanation
in terms of compactifications of buildings.
P. Bonvin was kind enough to write an appendix to this paper.
He wrote down the proof of a general commensurator super-rigidity 
which is used in Sect. 3.

\pec

Let us finally state a convention for group actions.
If a group $G$ acts on a set $X$, the (pointwise) stabilizer of a
subset $Y \subset X$ is called its {\it fixator~} and is denoted by 
${\rm Fix}_G(Y)$.
The classical (global) stabilizer is denoted by ${\rm Stab}_G(Y)$.
When $Y$ is a facet of a building $X$ and $G$ is a type preserving 
group of automorphisms, the
two notions coincide.
The notation $G \!\! \mid_Y$ refers to the group obtained from $G$ by 
factoring out
the kernel of the action on a $G$-stable subset $Y \subset X$.

\pec

I would like to thank M. Bourdon, M. Burger and Sh. Mozes for their 
constant interest in this approach to
Kac-Moody groups, as well as the audiences of talks at ETH Z\"urich 
(Advances in Rigidity Theory, June
2002, M. Burger and A. Katok organizers)  and at CIRM, Luminy 
(Geometry and Dynamics of Groups,
July 2002, M. Bourdon and L. Potyagailo organizers) for motivating 
questions and hints.
At last I am glad to express my deep gratitude to G. Rousseau; his 
careful reading of a previous version of
this paper enabled me to remove many inaccuracies and mistakes.

\vskip 15mm

\mysection{Structure theorem}
\label{s - Structure theorem}

We review some properties of topological Kac-Moody groups, keeping in 
mind the analogy with semisimple groups
over local fields of positive characteristic.
We quote results of combinatorial nature, to be used to prove the 
topological simplicity theorem, and we prove some results on
pro-$p$ Sylow subgroups.
Results on lattices are also shortly recalled.

\pec

\mysubsection{Topological Kac-Moody groups.---}
\label{ss - Topological KM groups}
Topological Kac-Moody groups were introduced in \cite[1.B]{RR02}.
Their geometric definition uses the buildings naturally associated to 
the countable Kac-Moody groups defined over finite fields.

\pec

\mysubsubsection{}
\label{sss - countable Tits} 
According to J. Tits \cite[3.6]{Tit87}, a {\it split Kac-Moody 
group~} is defined by generators
and Steinberg relations once a Kac-Moody root datum and a groundfield 
${\bf K}$ are given.
A {\it generalized Cartan matrix~} \cite{Kac90} is the main 
ingredient of a Kac-Moody root datum, the other part determining the
maximal split torus of the group.
More precisely, what J. Tits defines by generators and relations is a 
functor on rings.
When the generalized Cartan matrix is a Cartan matrix in the usual 
sense (which we call {\it of finite type}), the functor coincides
over fields with the functor of points of a Chevalley-Demazure group scheme.

\pec

An {\it almost split Kac-Moody group~} is the group of fixed points 
of a Galois action on a
split group \cite[\S 11]{Rem99}.   Let $\Lambda$ be an almost split 
Kac-Moody group.
Any such group satisfies the axioms of a {\it twin root 
datum~} \cite[\S 12]{Rem99}.
This is a group combinatorics sharply refining the structure of a 
$BN$-pair, and the associated geometry is a pair of twinned
buildings, conventionnally one for each sign \cite{Tit92}.
We denote by $X$ (resp. $X_-$) the positive (resp. negative) building 
of $\Lambda$.
The group $\Lambda$ acts diagonally on $X \times X_-$ in a natural way.
We choose  a pair of opposite chambers $R \subset X$ and $R_- \subset 
X_-$, which defines a pair of opposite
apartments $A \subset X$ and $A_- \subset X_-$  \cite[I.2]{Abr97}.
We call $R$ the {\it standard positive chamber~} and $A$ the  {\it 
standard positive apartment}.

\pec

{\sc Examples}.---~
1) The group ${\rm SL}_n({\bf K}[t,t^{-1}])$ is a Kac-Moody group of 
affine type.
The generalized Cartan matrix is $\widetilde A_{n-1}$ and the 
groundfield is ${\bf K}$.
More generally, values of Chevalley schemes over rings of Laurent 
polynomials are Kac-Moody
of affine type \cite[I.1 example 3]{Abr97}.

2) The group ${\rm SU}_3({\bf F}_q [t,t^{-1}])$ is an almost split 
Kac-Moody group of rank one.
The so-obtained geometry is a  semi-homogeneous twin tree of 
valencies $1+q$ and $1+q^3$
\cite[3.5]{Rem00}.

\pec

{\sc Definition}.---~\it
{\rm (i)~} We call $\Gamma$ the fixator of the negative chamber 
$R_-$, i.e. $\Gamma:={\rm Fix}_\Lambda(R_-)$.

{\rm (ii)~} We call $\Omega$ the fixator of the standard positive 
chamber $R$, i.e.
$\Omega:={\rm Fix}_\Lambda(R)$.
\rm\pec

{\sc Remark}.---~
In the above quoted references, the groups $\Lambda$, $\Gamma$ and 
$\Omega$ are denoted by
$G$, $B_-$ and $B$, respectively.
The reason is that the group $\Omega$ (resp. $\Gamma$) is the Borel 
subgroup of the
positive (resp. negative) Tits system of the twin root datum of 
$\Lambda$ \cite{Tit92},
\cite[\S 7]{Rem99}.

\pec

{\sc Example}.---~For $\Lambda={\rm SL}_n({\bf K}[t,t^{-1}])$, a 
natural choice of $R$ and $R_-$ defines
the group $\Omega$ as the subgroup of ${\rm SL}_n({\bf K}[t])$ which 
reduces to upper
triangular matrices modulo $t$, and $\Gamma$  is the subgroup of
${\rm SL}_n({\bf K}[t^{-1}])$ which reduces to lower triangular 
matrices modulo $t^{-1}$.

\pec

In this article, we are not interested in Tits group functors whose 
generalized Cartan matrices are of finite
type, whose values over finite fields are finite groups of Lie type.
Kac-Moody groups of affine type are seen as a guideline to generalize 
of classical results about
algebraic and arithmetic groups.
We are mainly interested in Kac-Moody groups which do not admit any 
natural matrix interpretation, and we want to understand to
what extent these new groups can be compared to the obviously linear ones.

\pec

\mysubsubsection{}
\label{sss - top KM} 
We now define topological groups -- see \cite[1.B]{RR02} for a wider framework.

\pec

{\sc Assumption}.---~\it
In this article, $\Lambda$ is a countable Kac-Moody group over the finite
field ${\bf F}_q$ of characteristic $p$ with $q$ elements, i.e. 
$\Lambda$ is the group
of rational points of an almost split Kac-Moody group with infinite 
Weyl group $W$.

\rm\pec

The kernel of the $\Lambda$-action on $X$ is the finite center
$Z(\Lambda)$, and the group $\Lambda/Z(\Lambda)$ still  enjoys the 
structure of twin root
datum \cite[Lemma 1.B.1]{RR02}.  Another consequence of the 
finiteness of the groundfield is
that the full automorphism group ${\rm Aut}(X)$ is  an uncountable 
totally disconnected
locally compact group.

\pec

{\sc Definition}.---~\it
{\rm (i)~} We call {\rm topological Kac-Moody group (associated to 
$\Lambda$)~} the closure in ${\rm Aut}(X)$ of the group
$\Lambda/Z(\Lambda)$.
We denote it by $\overline \Lambda$.

{\rm (ii)~}  We call {\rm parahoric subgroup (associated to $F$)~} 
the fixator in $\overline \Lambda$ of a given facet $F$.
We denote it by $\overline\Lambda_F$.
When the facet is a chamber, we call the corresponding group an {\rm 
Iwahori subgroup}.
\rm\pec

{\sc Remark}.---~ The group $\overline\Lambda$ is so to speak a 
completion of the group $\Lambda$.
Recall that the $\Lambda$-action on the single building $X$ is not 
discrete since the stabilizers are parabolic subgroups with
infinite \og unipotent radical\fg \cite[Theorem 6.2.2]{Rem99}.

\pec

{\sc Example}.---~For $\Lambda={\rm SL}_n({\bf K}[t,t^{-1}])$, $X$ is 
the Bruhat-Tits building of
${\rm SL}_n \bigl( {\bf F}_q(\!(t)\!) \bigr)$,  that is a Euclidean 
building of type $\widetilde A_{n-1}$.
If $\mu_n({\bf F}_q)$ denotes the $n$-th roots of unity in ${\bf 
F}_q$, the image $\Lambda/Z(\Lambda)$
of ${\rm SL}_n({\bf F}_q[t,t^{-1}])$ under the action on $X$ is ${\rm 
SL}_n({\bf F}_q[t,t^{-1}])/\mu_n({\bf F}_q)$ and the
completion $\overline\Lambda$ is
${\rm PSL}_n \bigl( {\bf F}_q(\!(t)\!) \bigr)={\rm SL}_n \bigl( {\bf 
F}_q(\!(t)\!) \bigr) /\mu_n({\bf F}_q)$.

\pec

\mysubsection{Refined Tits system and virtual pro-$p$-ness of 
parahoric subgroups.---}
\label{ss - Tits system}
The reference for this subsection is \cite[1.C]{RR02}.

\pec

\mysubsubsection{}
\label{sss - refined TS} 
Let us state  \cite[Theorem 1.C.2]{RR02} showing that topological 
Kac-Moody groups generalize
semisimple groups over local fields of positive characteristic.

\pec

{\sc Theorem}.---~\it
Let $\Lambda$ be an almost split Kac-Moody group over ${\bf F}_q$ and 
$\overline\Lambda$ be
its associated topological group.
Let $R \subset A$ be the standard chamber and apartment of the 
building $X$ associated to $\Lambda$.
We denote by ${\cal B}$ the standard Iwahori subgroup $\overline\Lambda_R$
and by $W_R$ the Coxeter group associated to $A$, generated by 
reflections along the panels of $R$.

\pec

{\rm (i)~} The topological Kac-Moody group $\overline\Lambda$ enjoys 
the structure of a refined
Tits system with abstract Borel subgroup ${\cal B}$ and Weyl group 
$W_R$, which is also
the  Weyl group of $\Lambda$.

{\rm (ii)~} Any spherical facet-fixator $\overline\Lambda_F$ is a semidirect 
product $M_F \ltimes
\widehat{U}_F$ where $M_F$ is a finite reductive group of Lie type 
over ${\bf F}_q$ and $\widehat{U}_F$ is a
pro-$p$ group.
In particular, any parahoric subgroup is virtually pro-$p$.
\qed\rm\pec

{\sc Remarks}.---~ 1) The group $\Lambda$ is {\it strongly transitive~} on the
building $X$, i.e. transitive on the pairs of chambers at given 
combinatorial (or $W$-) distance from one another
\cite[\S 5]{Ron89}.  This implies that $\overline\Lambda$ is strongly 
transitive on $X$ too, and
that $X$ is also the building associated to the above Tits system of 
$\overline\Lambda$.

2) Refined Tits system were defined by V. Kac and D. Peterson \cite{KP84}.
For twin root data, there are two relevant standard Borel subgroups playing
symmetric roles.
For refined Tits systems, only one conjugacy class of Borel subgroups 
is introduced.
The latter structure is abstractly implied by the former one 
\cite[1.6]{Rem99}, but it
applies to a strictly wider class of groups, e.g. ${\rm SL}_n \bigl( 
{\bf F}_q(\!(t)\!) \bigr)$
satisfies the axioms of a refined Tits system while it doesn't admit 
a twin root datum
structure.

\pec

Let us also briefly mention how the groups $M_F$ and $\widehat{U}_F$ 
are defined geometrically.
We first note that in the ${\rm CAT}(0)$ realization of buildings 
only spherical facets appear.
We denote by
${\rm St}(F)$ the {\it star~} of a facet $F$, that is the set of 
chambers whose closure contains
$F$.   Theorem 6.2.2 of \cite{Rem99} applies and we have a Levi decomposition
$\Lambda_F := {\rm Stab}_\Lambda(F) = M_F \ltimes U_F$, where $M_F$ 
is a Kac-Moody subgroup for a Cartan
submatrix of finite type and $U_F$ fixes pointwise ${\rm St}(F)$.
The group $\widehat{U}_F$ is the closure of $U_F$ in 
$\overline\Lambda$, hence it fixes ${\rm St}(F)$ too.
Moreover by  [loc. cit., Proposition 6.2.3], ${\rm St}(F)$ is a 
geometric realization of the
finite building attached  to $M_F$, and the action by $M_F$ is the 
standard one.
Therefore the image of the surjective homomorphism $\pi_F : 
\overline\Lambda_F \to M_F$ associated to
$\overline\Lambda_F=M_F \ltimes \widehat{U}_F$ gives the local action 
of $\overline\Lambda_F$
on ${\rm St}(F)$.

\pec

These facts are analogues of classical results in Bruhat-Tits theory 
\cite[Proposition
5.1.32]{BrT84}.  Namely, any facet in the Bruhat-Tits building of a 
semsimple group ${\bf G}$ over
a local field $k$ defines an integral structure over the valuation 
ring ${\cal O}$ of $k$.
The reduction of the ${\cal O}$-structure modulo the maximal ideal is 
a semisimple group over the residue
field, whose building is the star of the facet.
The integral points of the ${\cal O}$-structure act on it via the 
natural action of  the reduction.
The splitting $\overline\Lambda_F=M_F \ltimes \widehat{U}_F$ of a 
parahoric subgroup as a
semidirect product is specific to the case of valuated fields in 
equal characteristic, and in
the case of locally compact fields this only occurs in characteristic $p$.

\pec

{\sc Example}.---~
Let $v$ be a vertex in the Bruhat-Tits building of ${\rm SL}_n \bigl( 
{\bf F}_q(\!(t)\!) \bigr)$.
Then its fixator is isomorphic to ${\rm SL}_n({\bf F}_q[[t]])$ and 
its star is isomorphic to the building of ${\rm SL}_n({\bf F}_q)$.
The subgroup $\widehat{U}_v$ is the first congruence subgroup of 
${\rm SL}_n({\bf F}_q[[t]])$, i.e. the matrices reducing
to the identity modulo $t$.
The above reduction corresponds concretely to moding out by 
$\widehat{U}_v$, and the Iwahori subgroups
fixing a chamber in ${\rm St}(v)$ are the subgroups reducing to a 
Borel subgroup of ${\rm SL}_n({\bf F}_q)$, e.g. reducing to
the upper triangular matrices.

\pec

In our case, the Bruhat decomposition and the rule to multiply double classes
\cite[IV.2]{Bou81}  have topological consequences.

\pec

{\sc Corollary}.---~\it
A topological Kac-Moody group is compactly generated.
\rm\pec

{\it Proof}.---~ 
Let $\overline\Lambda$ be such a group and let ${\cal 
B}=\overline\Lambda_R$ be the standard Iwahori subgroup.
Then $\overline\Lambda$ is generated by ${\cal B}$ and by the compact 
double classes ${\cal B}s{\cal B}$, when
$s$ runs over the finite set of reflections along the panels of $R$.
\qed\pec

\mysubsubsection{}
\label{sss - Sylow} 
Pro-$p$ Sylow subgroups of totally disconnected groups are defined in 
\cite[I.1.4]{Ser94}  for
instance.
The following proposition is suggested by classical results on 
pro-$p$ subgroups of
non-Archimedean Lie groups, e.g. \cite[Theorem 3.10]{PR94}.

\pec

{\sc Proposition}.---~\it
{\rm (i)~} Given any chamber $R$, the group $\widehat{U}_R$ of the 
decomposition
$\overline\Lambda_R=M_R \ltimes \widehat{U}_R$ is the unique pro-$p$ 
Sylow subgroup of the Iwahori subgroup
$\overline\Lambda_R$.

{\rm (ii)~} Let $K$ be a pro-$p$ subgroup of $\overline\Lambda$.
Then there is a chamber $R$ such that $K$ lies in the pro-$p$ Sylow 
subgroup $\widehat{U}_R$ of
$\overline\Lambda_R$.

{\rm (iii)~} The pro-$p$ Sylow subgroups of $\overline\Lambda$ are 
precisely the subgroups $\widehat{U}_R$
when $R$ ranges over the chambers of the building $X$; they are all conjugate.

{\rm (iv)~} The Iwahori subgroups of $\overline\Lambda$ are 
intrinsically characterized as the normalizers of the pro-$p$ Sylow
subgroups of $\overline\Lambda$.
\rm\pec

{\it Proof}.---~
(i). By quasi-splitness of an almost split Kac-Moody group $\Lambda$ 
over ${\bf F}_q$
\cite[13.2]{Rem99}, the Levi factor $M_R$ of a chamber fixator in 
$\Lambda$ is the ${\bf
F}_q$-points of a torus.  Therefore its order is prime to $p$ and we 
conclude by the
decomposition
$\overline\Lambda_R=M_R \ltimes \widehat{U}_R$.

\pec

(ii). Let $K$ be a pro-$p$ subgroup of $\overline\Lambda$.
Since it is compact, it fixes a spherical facet $F$ \cite[4.6]{Rem99} 
and by Theorem \ref{sss -
refined TS} (ii)  we can write $K < \overline\Lambda_F = M_F \ltimes 
\widehat{U}_F$.
Let us look at the local action $\pi_F :  \overline\Lambda_F \to M_F$ 
(\ref{sss - refined TS}).
By the Bruhat decomposition in split $BN$-pairs, the $p$-Sylow 
subgroups of finite reductive
groups of Lie type are the unipotent radicals of the Borel subgroups
\cite[B Corollary 3.5]{Bor70}, so the $p$-subgroup $\pi_F(K)$ of 
$M_F$ is contained in some Borel
subgroup of $M_F$,  hence fixes a chamber $R$ in ${\rm St}(F)$. 
Since $\widehat{U}_F$ fixes
pointwise ${\rm St}(F)$, the whole subgroup $K$ fixes $R$ and we 
conclude by (i).

\pec

(iii). The first assertion follows immediately from (ii), and the 
second one follows from the transitivity of $\Lambda$ on the
chambers of $X$.

\pec

(iv). According to (iii), it is enough to show that we have ${\cal B} 
= N_{\overline\Lambda}(\widehat{U}_R)$
for the Iwahori subgroup ${\cal B}=\overline\Lambda_R$ fixing the 
standard positive chamber $R$.
By \ref{sss - refined TS} (ii), we know that ${\cal B} < 
N_{\overline\Lambda}(\widehat{U}_R)$.
By \cite[IV.2.5]{Bou81} the normalizer 
$N_{\overline\Lambda}(\widehat{U}_R)$ is an abstract
parabolic subgroup of the Tits system of \ref{sss - refined TS} (i) 
with abstract Borel subgroup
${\cal B}$.  If $N_{\overline\Lambda}(\widehat{U}_R)$ were bigger 
than ${\cal B}$, it would
contain a reflection conjugating a positive root group to a negative one.
But this is in contradiction with axiom (RT3) of refined Tits systems 
\cite{KP84}, so we have
${\cal B}=N_{\overline\Lambda}(\widehat{U}_R)$.
\qed
\pec

{\sc Remark}.---~ 
M. Sapir pointed out to us that knowing whether the pro-$p$ Sylow 
subgroups of some $\overline\Lambda$ are new is an interesting
question too.
Note that it is a very hard problem (solved by E. Zelmanov) to show 
that the free pro-$p$ group is
non-linear in low dimension.

\pec

\mysubsection{Lattices and generalized arithmeticity.---}
\label{ss - Lattices}
We briefly discuss existence and generalized arithmeticity of 
lattices in topological Kac-Moody groups.

\pec

\mysubsubsection{}
\label{sss - KM lattices} 
We keep the almost split Kac-Moody group $\Lambda$ over ${\bf F}_q$.
The Weyl group $W$ is infinite, and we denote by $W(t)$ its growth 
series $\sum_{w \in W}  t^{\ell(w)}$.

\pec

{\sc Theorem}.---~\it
Assume that $W({1/q}) < \infty$.
Then $\Lambda$ is a lattice of $X \times X_-$ for its diagonal 
action, and for any
point $x_- \! \in \! X_-$ the subgroup $\Lambda_{x_-}={\rm 
Fix}_\Lambda(x_-)$ is a lattice of $X$.
These lattices are not uniform.
\qed\rm\pec

This result was independently proved in \cite{CG99} and in \cite{Rem98}.
In the case of ${\rm SL}_n({\bf F}_q[t,t^{-1}])$, the lattices of the 
form $\Lambda_{x_-}$ are all commensurable to the
arithmetic lattice ${\rm SL}_n({\bf F}_q[t^{-1}])$ in ${\rm 
SL}_n({\bf F}_q (\!( t )\! )$.
Recall that by Margulis commensurator criterion \cite[Theorem 6.2.5]{Zim84} ,
a lattice in a semisimple Lie group $G$ is arithmetic if and only if 
its commensurator is dense in $G$.
Taking this characterization as a definition for general situations,
\cite[Lemma 1.B.3 (ii)]{RR02} says that the groups $\Lambda_{x_-}$ are
arithmetic lattices of $\overline\Lambda$ by the very definition of 
this topological group.
Here is the statement, whose proof is based on refined Tits sytem arguments.
\pec

{\sc Lemma}.---~\it
For any $x_- \! \in \! X_-$, we have:
$\Lambda<{\rm Comm}_{\overline\Lambda}(\Lambda_{x_-})$.
\qed\rm\pec

{\sc Remark}.---~ 
Some lattices may have big enough commensurators to be arithmetic in 
${\rm Aut}(X)$, meaning that the
commensurators are dense in ${\rm Aut}(X)$.
This is the case for the Nagao lattice ${\rm SL}_2({\bf 
F}_q[t^{-1}])$ in the full automorphism group of the $q+1$-regular 
tree
\cite{Moz00}.
The proof can be formalized  and extended to exotic trees admitting a 
Moufang twinning
\cite{AR02}.

\pec

\mysubsubsection{}
\label{sss - uniform lattices} 
A way to produce lattices in automorphism groups of cell-complexes is 
to take fundamental groups of complexes of
groups \cite[III.${\cal C}$]{BH99}, the point being then to recognize 
the covering space.
A positive result is the following -- see \cite{Bou00}: let $R$ be a 
regular right-angled polygon in the hyperbolic plane ${\bf H}^2$
and $\underline q:= \{ q_i \}_{1 \leq i \leq r}$ be a sequence of 
integers $\geq 2$.
(When $\underline q$ is constant, we replace $\underline q$ by its value $q$.)
Then there exists a unique {\it right-angled Fuchsian building~} 
$I_{r,1+\underline q}$ with apartments isomorphic to the tiling
of ${\bf H}^2$ by $R$, and such that the link at any vertex of type 
$\{ i; i+1 \}$ is the complete bipartite graph of parameters
$(1+q_i,1+q_{i+1})$.
The so-obtained lattices are uniform and abstractly defined by:
$\Gamma_{r,1+\underline q}=\langle \, \{ \gamma_i \}_{i \in {\bf 
Z}/r} \,  : \, \gamma_i^{q_i+1}=1$ and
$[\gamma_i,\gamma_{i+1}]=1 \, \rangle$.
This uniqueness is a key argument to prove \cite[Proposition 5.C]{RR02}:

\pec

{\sc Proposition}.---~\it
For any prime power $q$, there exists a non-uniquely defined 
Kac-Moody group $\Lambda$ over ${\bf
F}_q$ whose building is $I_{r,1+q}$.
Moreover we can choose $\Lambda$  such that its natural image in 
${\rm Aut}(I_{r,1+q})$ contains the uniform lattice
$\Gamma_{r,1+q}$.
\qed\rm\pec

This result says that the buildings $I_{r,1+q}$ are relevant to both 
Kac-Moody theory, and generalized hyperbolic geometry since
they carry a natural ${\rm CAT}(-1)$-metric.
They provide the most well-understood infinite family of exotic 
Kac-Moody buildings (indexed by $r\geq 5$ when $q$ is fixed).
The corresponding countable groups $\Lambda$ are actually typical 
groups to which our non-linearity result applies (\ref{ss -
non-linear}).
We study these buildings more carefully in \ref{ss - Fuchsian}.
Finally, combining  \cite[Theorem 4.6]{Rem01} and \cite{Bou97} leads 
to a somewhat surprising
situation, with coexistence of lattices with sharply different properties
\cite[Corollary 5.C]{RR02}.

\pec

{\sc Corollary}.---~\it
Whenever $q$ is large enough and $r$ is even and $> \!\! > q$, the 
topological group $\overline\Lambda$ associated to the above
$\Lambda$ contains both uniform Gromov-hyperbolic lattices  embedding 
convex-cocompactly into real hyperbolic spaces, and
non-uniform Kac-Moody lattices which are not linear in characteristic 
$\neq p$, containing infinite groups
of exponent dividing $p^2$.
\qed\rm\pec

{\sc Remark}.---~
Determining the commensurator of the uniform lattice 
$\Gamma_{r,1+\underline q}$, hence
deciding its arithmeticity in ${\rm Aut}(I_{r,1+\underline q})$, is 
an open question.

\gec

\mysection{Topological simplicity theorem}
\label{s - Topological simplicity theorem}

We prove that topological Kac-Moody groups as previously defined are 
direct products of topologically simple groups.
In view of existence of congruence subgroups in the affine case, the 
proof is expected not to work for countable Kac-Moody groups.
A failing argument and some open problems are discussed, and some 
results on Levi factors and homomorphisms to non-Archimedean
groups are proved.

\pec

\mysubsection{Topological simplicity.---}
\label{ss - Topological simplicity}
Here is a further argument supporting the analogy with semisimple 
algebraic groups over local
fields of positive characteristic.

\pec

\mysubsubsection{}
\label{sss - top simple}
As for simplicity of classical groups, we need to assume the 
groundfield large enough, because in our proof we
need simplicity of some rank-one finite groups of Lie type.

\pec

{\sc Theorem}.---~\it
Let $\Lambda$ be a countable Kac-Moody group which is almost split 
over the finite field ${\bf F}_q$, with
$q>4$.  We assume that $\Lambda$ is generated by its root groups, 
e.g. simply connected.

\pec

{\rm (i)~} If the Dynkin diagram of $\Lambda$ is connected, the 
associated topological Kac-Moody group $\overline\Lambda$ is
topologically simple.

{\rm (ii)~} For an arbitrary Dynkin diagram, the group 
$\overline\Lambda$ is a direct product of topologically simple 
groups, each
factor being the topological Kac-Moody group associated to a 
connected component of the Dynkin diagram.
\rm\pec

{\it Proof}.---~
(ii). The building of a Kac-Moody group $\Lambda$ is defined as a gluing
$\displaystyle X:={\Lambda \times A \over \sim}$, where $A$ is the model for an
apartment, i.e. the ${\rm CAT}(0)$ realization of the Coxeter complex 
of the Weyl group $W$.
The $\Lambda$-action comes from factorizing the map
$\Lambda \times \Lambda \times A \to \Lambda \times A$ which sends 
$(\lambda',\lambda,x)$ to
$(\lambda'\lambda,x)$ \cite[\S 4]{Rem99}.
The rule by which the Coxeter diagram of $W$ is deduced from the Dynkin diagram
\cite[3.1]{Tit87}  implies that irreducible factors correspond to 
connected components of the
diagram.  The model $A$ is then the direct product of the models for 
the Coxeter complexes of the
irreducible factors of $W$.
By the defining relations of $\Lambda$ \cite[3.6]{Tit87},
a root group indexed by a root in the subsystem associated to a given 
connected component of the Dynkin
diagram centralizes a root group arising from another connected component.
By the relationship between buildings and Tits systems \cite[\S 5]{Ron89},
the Kac-Moody subgroup defined by  a given connected component acts
trivially on a factor of the building $X$ arising from any other 
connected component.
Therefore proving (ii) is reduced to proving (i).

\pec

(i). Let ${\cal B}=\overline \Lambda_R$ be the standard Iwahori subgroup.
By Theorem \ref{sss - refined TS} (i), it is the Borel subgroup of a 
Tits system with the
same Coxeter system as the one for $\Lambda$.
By the Kac-Moody analogue of Lang's theorem \cite[13.2.2]{Rem99}, the 
Levi factor $M_R$ in ${\cal
B}=M_R
\ltimes \widehat{U}_R$ is a maximally split maximal torus $T$ of 
$\Lambda$, i.e. the rational
points of a finite
${\bf F}_q$-torus.

\pec

Let $I$ be the indexing set of the simple roots of $\Lambda$ and let 
$G_i$ be the standard semisimple Levi factor
of type $i \! \in \! I$.
The group $G_i$ is a finite almost simple group of Lie type generated 
by the root groups
attached to the simple root $a_i$ and its  opposite \cite[6.2]{Rem99}.
By our assumption they generate $\Lambda$ as
an abstract group, hence $\overline\Lambda$ as a topological group.
Note that $G_i$ has no non-trivial abelian quotient, and it has no 
quotient  isomorphic to a
$p$-group either -- see Remark 2 below.

\pec

We isolate the remaining arguments in a lemma also used to prove 
\ref{sss - Levi factors} (iv).

\pec

{\sc Lemma}.---~\it
Let $G$ be a topological group acting continuously and strongly 
transitively by type-preserving automorphisms on a building $X$
with irreducible Weyl group.
We denote by ${\cal B}$ a chamber fixator and we assume that it is an 
abelian-by-pro-$p$ extension.
We also assume that $G$ is topologically generated by finite 
subgroups admitting no non-trivial quotient
isomorphic to any abelian or $p$-group.
Then any proper closed normal subgroup of $G$ acts trivially on $X$.
\rm\pec

This proves the theorem by setting $G:=\overline\Lambda$, ${\cal 
B}:=\overline\Lambda_R=T
\ltimes \widehat{U}_R$ and by choosing $\{ G_i \}_{i \in I}$ as 
generating family of
subgroups, since by definiton $\overline\Lambda$ acts faithfully on $X$.
\qed\pec

Let us now prove the lemma.

\pec

{\it Proof}.---~
By \cite[\S 5]{Ron89}, $G$ admits an irreducible Tits system with 
${\cal B}$ as abstract Borel
subgroup.  Assume we are given a closed normal subgroup $H \triangleleft_c G$.
Then by \cite[IV.2.7]{Bou81}, we have either
$H.{\cal B}={\cal B}$ or $H.{\cal B}=G$.
The first case implies $H < {\cal B}$, and since $H$ is normal in $G$:

\pec

\centerline{$\displaystyle H < \bigcap_{g \in G} g{\rm Fix}_G(R)g^{-1}
={\rm Fix}_G \bigl( \bigcup_{g \in G} gR \bigr)
={\rm Fix}_G \bigl( \bigcup_{g \in G} g\overline R \bigr)= {\rm Fix}_G(X)$,}

\pec

because $G$ is type-preserving and transitive on the chambers of the 
building $X$.

\pec

 From now on, we assume that $H$ doesn't act trivially on $X$.
By the previous point, this implies that we have the identification 
of compact groups:
$G/H \simeq {\cal B}/({\cal B}\cap H)$.
We denote by $1 \to \widehat{U} \to {\cal B} \to T \to 1$ the 
extension of the assumption,
and we consider the homomorphism $\widehat{U} \to {\cal B}/({\cal 
B}\cap H)$ sending $u$ to
$u({\cal B}\cap H)$.
Its kernel is $\widehat{U} \cap H$, so we have an injection
$\widehat{U}/(\widehat{U}\cap H) < {\cal B}/({\cal B}\cap H)$, and
since $\widehat{U}$ is normal in ${\cal B}$, so is 
$\widehat{U}/(\widehat{U}\cap H)$ in
${\cal B}/({\cal B}\cap H)$.
Then we consider the composition of surjective homomorphisms
$\displaystyle {\cal B} \to {\cal B}/({\cal B}\cap H) \to
{{\cal B}/({\cal B}\cap H) \over \widehat{U}/(\widehat{U}\cap H)}$.
Its kernel contains $\widehat{U}$, which shows that
$\displaystyle {{\cal B}/({\cal B}\cap H) \over 
\widehat{U}/(\widehat{U}\cap H)}$ is a quotient
of the abelian group $T$.
Therefore if we denote by $\{ G_i \}_{i \in I}$ the generating family 
of subgroups of the assumptions,
the image of each $G_i$ under
$\displaystyle G \to G/H \simeq {\cal B}/({\cal B}\cap H)
\to {{\cal B}/({\cal B}\cap H) \over \widehat{U}/(\widehat{U}\cap 
H)}$ is trivial.

\pec

In other words, the image of $G_i$ under
$\displaystyle G \to G/H \simeq {\cal B}/({\cal B}\cap H)$
is a finite subgroup of the  group
$\displaystyle \widehat{U}/(\widehat{U}\cap H)$.
But the group $\widehat{U}/(\widehat{U}\cap H)$ is pro-$p$ since it 
is the quotient of a pro-$p$ group
by a closed normal subgroup \cite[1.11]{DdSMS99}.
Therefore, once again in view of the possible images of $G_i$, we 
must have: $G_i/(G_i \cap H) = \{ 1 \}$, that
is $G_i < H$.
Since the groups $G_i$ topologically generate $G$, we have:  $H=G$.
\qed
\pec

{\sc Example}.---~
In order to see which classical result is generalized here,
we take $\Lambda={\rm SL}_n({\bf F}_q[t,t^{-1}])$.
Let $\mu_n({\bf F}_q)$ be the  $n$-th roots of unity in
${\bf F}_q$.  The image of ${\rm SL}_n({\bf F}_q[t,t^{-1}])$ for the 
action on its positive
Euclidean building is
${\rm SL}_n({\bf F}_q[t,t^{-1}])/\mu_n({\bf F}_q)$.
Though ${\rm SL}_n({\bf F}_q[t,t^{-1}])$ does not come from a simply 
connected Kac-Moody root
datum, it is generated by standard rank-one Levi factors.
The lemma says that the closed normal subgroups of
${\rm SL}_n \bigl({\bf F}_q(\!( t )\!) \bigr)$ are central, i.e. in $\mu_n({\bf 
F}_q)$.

\pec

{\sc Remarks}.---~
1) An argument in the proof is that the ${\bf F}_q$-almost simple 
finite Levi factors of the panel
fixators generate the countable group $\Lambda$ and don't admit any 
quotient isomorphic to an abelian or a
$p$-group.
The same proof works if the panel fixators are replaced by a 
generating family of facet
fixators with the same property on quotients.
For ${\rm SL}_3({\bf F}_2(\!(t)\!))$ the Levi factors of the panels 
fixators are solvable
$\simeq {\rm SL}_2({\bf F}_2)$, but the proof works with the Levi 
factors for the three
standard vertices, because these groups are isomorphic to ${\rm 
SL}_3({\bf F}_2)$.

2) Assuming $q>4$ is a way to have this property on quotients for any 
facet fixator since
non-simplicity of rational points of adjoint ${\bf F}_q$-simple 
groups occurs only for
$q \leq 4$ \cite[11.1.2 and 14.4.1]{Car89}.
Indeed, let $G$ be the group generated by the $p$-Sylow subgroups of 
the rational points of an
almost simple ${\bf F}_q$-group with $q>4$, and let $H=G/K$ be a quotient.
By simplicity, $K/Z(G)$ equals $\{ 1 \}$ or $G/Z(G)$.
The first case implies that $H$ surjects onto the non-abelian simple 
group $G/Z(G)$.
In the second case we have $G=K\cdot Z(G)$, showing that $H$ is a 
quotient of $Z(G)$ of order
prime to $p$.
If $U$ denotes a $p$-Sylow subgroup of $G$ then $U/(U \cap K)$ is 
trivial in $H$, implying
that $K$ contains the $p$-Sylow subgroups of $G$, hence equals $G$.

\pec

\mysubsubsection{}
\label{sss - trivial image}
In \cite{Rem01} it is proved that the group $\Lambda$ cannot be 
linear in characteristic $\neq p$ (e.g. in characteristic 0),
meaning that no homomorphism $\Lambda \to {\rm GL}_N(k)$ hence no homomorphism
$\overline \Lambda \to {\rm GL}_N(k)$ can be injective.
Combined with the topological simplicity of $\overline\Lambda$, this leads to:

\pec

{\sc Corollary}.---~\it
Let $k$ be a local field of characteristic $\neq p$ and let ${\bf G}$ 
be a linear algebraic group over $k$.
Then any continuous homomorphism $\varphi: \overline \Lambda \to {\bf 
G}(k)$ has trivial image.
\qed\rm\pec

{\sc Remark}.---~ 
According to \cite[Theorem 4.E.2]{RR02}, there exist generalizations 
of Kac-Moody groups with
Weyl groups of arbitrarily large rank and with several groundfields.
The consequence of the latter point is the complete non-linearity of 
some of these groups  [loc.
cit., Theorem 5.B].  It is expectable that the corresponding 
topological groups are topologically
simple, but the pro-$p$-ness argument in the proof of Lemma \ref{sss 
- top simple} is not
available.  This would lead to the same corollary, yet without any 
restriction on the target
characteristic.

\pec

\mysubsection{Discussion of the proof for countable groups.---}
\label{ss - Discussion}
The proof of Lemma \ref{sss - top simple} is inspired by 
\cite[IV.2.7]{Bou81}, where an abstract
simplicity criterion is derived from properties of Tits systems.
We show why the proof doesn't work for countable Kac-Moody groups, 
which will lead us to natural questions about
the latter groups.

\pec

\mysubsubsection{}
\label{sss - discrete congruence} 
The basic idea is to require the irreducibility of a Tits system, and 
two further assumptions involving the
\og unipotent radical\fg $U$ of the abstract Borel subgroup.
The first one is a generation assumption on $U$ and commutator subgroups.
The second one -- property (R) of \cite[IV.2.7]{Bou81} -- is 
technical, and satisfied by
solvable or simple groups.  A consequence of this result is a 
conceptual proof of the abstract
simplicity of ${\rm PSL}_n$ over (large enough) fields.

\pec

{\sc Remark}.---~ 
In our situation, the Borel subgroup of the Tits system is ${\cal B}$ 
and the pro-$p$ group
$\widehat{U}_R \triangleleft_{\rm f.i.} \overline\Lambda_R$ is the 
\og abstract unipotent radical\fg.
Virtual pro-$p$-ness of the Iwahori subgroup ${\cal 
B}=\overline\Lambda_R$ replaces property (R).
We don't know if checking the assumptions of \cite[IV.2.7]{Bou81}  is 
doable; which would provide
an abstract simplicity theorem.

\pec

Back to $\Lambda={\rm SL}_n({\bf F}_q[t,t^{-1}])$, we see that we 
cannot apply the arguments of the above proof to it
since the conclusion is false: this arithmetic group is {\it 
residually finite}, which means
that the intersection
$\bigcap_{\Lambda' \triangleleft_{\rm f.i.} \Lambda} \Lambda'$
of its finite index normal subgroups is trivial.
A way to see this is to use {\it congruence subgroups~} $K({\goth 
a})$, where ${\goth a}$ is an ideal of
${\bf F}_q[t,t^{-1}]$.
The congruence subgroup of ${\rm SL}_n({\bf F}_q[t,t^{-1}])$ attached 
to ${\goth a}$ is
by definition the kernel of the reduction homomorphism:

\pec

\centerline{${\rm SL}_n({\bf F}_q[t,t^{-1}])
\to {\rm SL}_n({\bf F}_q[t,t^{-1}]/{\goth a})$.}

\pec

Congruence subgroups enable to see on this example what goes wrong in the above
proof when $\overline \Lambda$ is replaced by $\Lambda$.
A key point for $\overline \Lambda$ is:

\pec

\hskip 4mm $(\star)$ a finite rank-one simple Levi factor
cannot be mapped non-trivially into $\widehat{U}/(\widehat{U} \cap H)$.

\pec

This is proved by using the facts that the quotient of a pro-$p$ 
group by a normal closed subgroup is pro-$p$ and that
pro-$p$-ness  implies  strong  restrictions on closed, e.g. finite, subgroups.
Here the pro-$p$ group $\widehat{U}_R \triangleleft_{\rm f.i.} {\cal B}$ is:

\pec

\centerline{$\widehat{U}_R=
{\rm SL}_n \pmatrix{\xymatrix{1+t{{\bf F}_q[[t]]} \ar@{.}[dr] & {{\bf 
F}_q[[t]]}
\\ t{{\bf F}_q[[t]]} & 1+t{{\bf F}_q[[t]]}}}$.}

\pec

It is the closure of the subgroup:

\pec

\centerline{$U={\rm SL}_n \pmatrix{\xymatrix{1+t{{\bf F}_q[t]} 
\ar@{.}[dr] & {{\bf F}_q[t]}
\\ t{{\bf F}_q[t]} & 1+t{{\bf F}_q[t]}}}$,}

\pec

which replaces $\widehat{U}_R$ in the discrete case.
The ideal ${\goth a}=(t-1){\bf F}_q[t]$ provides a congruence 
subgroup $K({\goth a})$.
Moding out by $K({\goth a})$ amounts to applying the identification 
$t \sim 1$, for which the
quotient of $U$ is ${\rm SL}_n({\bf F}_q)$, hence contains several 
copies of ${\rm SL}_2({\bf F}_q)$.
This shows that $U$ contains normal subgroups $H$ such that $U/(U 
\cap H)$ contains finite rank-one simple groups of Lie type.
Thus the above quoted key fact $(\star)$ doesn't hold when replacing 
$\overline\Lambda$ (resp. $\widehat{U}_R$) by $\Lambda$
(resp. $U$).

\pec

\mysubsubsection{}
\label{sss - question discrete} 
In view of the existence of congruence subgroups e.g. in ${\rm 
SL}_n({\bf F}_q[t,t^{-1}])$, it is worth discussing for a general
$\Lambda$ a group-theoretic property that is opposite to the notion 
of simplicity.
The affine case thus leads to the following

\pec

{\sc Question}.---~\it
Which countable Kac-Moody groups are residually finite?
\rm\pec

{\sc Remark}.---~ 
This question was asked by M. Burger, who notes that taking fixators 
of combinatorial balls centered at a
chamber $R$ provides a family of finite index normal subgroups of 
$\Omega={\rm Fix}_{\Lambda}(R)$.
Since the intersection of these groups is trivial, $\Omega$ is 
residually finite, and
the same argument works for $\Gamma$ using the negative building.

\pec

According to  \cite[Proposition 1.7]{Kac90}, if a generalized Cartan 
matrix is indecomposable and
invertible, then the corresponding Kac-Moody algebra over ${\bf C}$ 
is simple, but
the relation between Kac-Moody algebras and groups is very loose.
Still, the fact that affine Kac-Moody algebras are very specific for 
simplicity suggests
the following

\pec

{\sc Question}.---~\it
Are there simple infinite countable Kac-Moody groups?
\rm\pec

{\sc Remark}.---~ 
According to Mal'cev theorem \cite[\S 3]{Mag69}, proving the non 
residual finiteness (and {\it a
fortiori~} the simplicity) of an infinite Kac-Moody group $\Lambda$ 
is a way to disprove any
linearity.
The previous subsection shows that the method of \ref{sss - top 
simple} cannot be applied to
countable Kac-Moody groups, but for non residual finiteness, another 
strategy could be to use
analogues of the criterion introduced  by M. Burger and Sh. Mozes in 
the context of products of
trees \cite[\S 2]{BM00}.

\pec

\mysubsection{Closures of Levi factors and maps to non-Archimedean groups.---}
\label{ss - Levi and maps}
We investigate the structure of closures of Levi factors.
This enables us to prove a result on continuous homomorphisms from 
topological Kac-Moody groups to Lie groups over
non-Archimedean fields.

\pec

\mysubsubsection{}
\label{sss - Levi factors}
We keep our Kac-Moody group $\Lambda$ and the inclusion of the 
chamber $R$ in the apartment $A$.
We denote by $\{ a_i \}_{i \in I}$ the finite family of simple roots, 
which we see as half-spaces of $A$ whose intersection equals
$R$.
Let us now choose a subset $J$ of $I$.
We can thus introduce the {\it standard parabolic subgroup 
$\Lambda_J$~} in $\Lambda$, which
is the union
$\Lambda_J=\bigsqcup_{w \in W_J} \Omega w \Omega$ indexed by the 
Coxeter group generated by the reflections along the walls
$\partial a_i$ for $i \! \in \! J$.
By \cite[Theorem 6.2.2]{Rem99}, the group $\Lambda_J$ admits a Levi 
decomposition $\Lambda_J = M_J
\ltimes U_J$, where the Levi factor $M_J$ is the Kac-Moody subgroup 
generated by the maximal
split torus $T$ and the root groups indexed by the roots
$w.a_i$ for $w \! \in \! W_J$ and $i \! \in \! J$.

\pec

{\sc Definition}.---~\it
{\rm (i)~} We denote by $\overline M_J$ the closure of $M_J$ in 
$\overline \Lambda$, and by
$\overline G_J$ the topological group generated by the root groups 
$U_{\pm a_i}$ when $i$ ranges over $J$.

{\rm (ii)~} The intersection of roots $R_J:=\bigcap_{i \in J} a_i$ is 
called the {\rm inessential chamber~} of $\overline M_J$ in
$X$.

{\rm (iii)~} We denote by ${\cal B}_J$ the stabilizer of $R_J$ in 
$\overline M_J$.

{\rm (iv)~} The union of closures $X_J:=\bigcup_{g \in M_J} 
g.\overline R_J$ is called the {\rm inessential building~} of
$\overline M_J$ in $X$.
\rm\pec

What supports the terminology is the following

\pec

{\sc Proposition}.---~\it
{\rm (i)~} The closure group $\overline M_J$ admits a natural refined 
Tits system with Weyl group $W_J$ and abstract Borel
subgroup ${\cal B}_J$.

{\rm (ii)~} The subgroup ${\cal B}_J$ fixes pointwise $\overline R_J$ 
and admits a decomposition
${\cal B}_J = T \ltimes \widehat{U_J}$, where $T$ is the maximally 
split maximal
${\bf F}_q$-torus attached to $A$ and $\widehat{U_J}$ is pro-$p$.

{\rm (iii)~} The space $X_J$ is a geometric realization of the 
building of $\overline M_J$ arising from the above Tits system
structure.

{\rm (iv)~} The group $\overline G_J$ is of finite index in 
$\overline M_J$ and when
$q>4$ it admits a direct product decomposition into topologically 
simple factors.
\rm\pec

{\it Proof}.---~ 
(ii). By \cite[6.2]{Rem99}, the group $M_J$ has a twin root datum 
structure with positive Borel
subgroup $\Omega_J:=\Omega\cap M_J$, negative Borel subgroup 
$\Gamma_J:=\Gamma\cap M_J$ and
Weyl group $W_J$.  Moreover the group $\Omega_J$ is generated by $T$ 
and the root groups $U_a$
where $a$ is a positive root of the form
$w.a_i$ with $i \! \in \! J$ and $w \! \in \! W_J$.
By the proof of Corollary 1 in \cite[5.7]{MP95}, such a root contains 
$R_J$, so the corresponding
group $U_a$ fixes
$\overline R_J$ pointwise.
Since $T$ fixes pointwise the whole apartment $A$ and ${\cal 
B}_J=\overline{\Omega_J}$, we proved the first assertion.
We have $T < \Omega_J$ and ${\cal B}_J<{\cal B}$ since ${\cal B}_J$ 
fixes $R$: this proves the second point by setting
$\widehat{U_J}:={\cal B}_J \cap \widehat{U}_R$.

\pec

The Borel subgroup $\Omega_J$ of the positive Tits system of $M_J$ is 
the fixator of $R_J$, so $X_J$ is a geometric
realization of the corresponding building and $\overline M_J$ is 
strongly transitive on $X_J$.
This proves (iii) and the fact that $\overline M_J$ admits a natural 
Tits system with Weyl group $W_J$ and
abstract Borel subgroup ${\cal B}_J$.
We are now in position to apply the same arguments as in 
\cite[1.C]{RR02}  to prove that we
actually have a refined Tits system,  which proves (i).

\pec

(iv). The first assertion follows from the fact that we have 
$\overline M_J=\overline G_J \cdot T$.
The same argument as for Theorem \ref{sss - top simple} (ii) implies 
that we are reduced to considering the groups corresponding
to  the connected components of the Dynkin subdiagram obtained from 
the one of $\Lambda$ by removing the vertices of type
$\not\in J$ (and the edges containing one of them).
Such groups commute with one another, and they are topologically 
simple by Lemma \ref{sss - top simple}.
\qed\pec

{\sc Remarks}.---~ 
1) Any proper submatrix of an irreducible affine  generalized Cartan 
matrix is of finite type
\cite[Theorem 4.8]{Kac90}.  Consequently, the Levi factors of the 
parabolic subgroups of
$\Lambda$ are all finite groups of Lie type in this case.

2) Conversely, if the submatrix of type $J$ of the generalized Cartan 
matrix of $\Lambda$ is not spherical (i.e. of finite type), then
$\overline M_J$ contains an infinite topologically simple subgroup.
Such a  group comes from a non spherical connected component of the 
Dynkin subdiagram given by $J$.

\pec

\mysubsubsection{}
\label{sss - maps} 
The virtual topological simplicity of closures of Levi factors 
enables to prove the following
result about actions of Kac-Moody groups on Euclidean buildings.

\pec

{\sc Proposition}.---~\it
Let $\Lambda$ be a countable Kac-Moody group which is almost split 
over the finite field
${\bf F}_q$ with $q>4$, which is generated by its root groups and 
which has connected Dynkin
diagram.
We assume we are given a continuous group homomorphism
$\mu : \overline\Lambda \to {\bf G}(k)$, where ${\bf G}$ is a 
semisimple group defined over a
non-Archimedean local field $k$.
We denote by $\Delta$ the Bruhat-Tits building of ${\bf G}$ over $k$, 
and for each point
$x$ in the Kac-Moody building $X$, we denote by $\overline\Lambda_x$ 
the fixator
${\rm Fix}_{\overline\Lambda}(x)$.

\pec

{\rm (i)~} For any point $x$ in the Kac-Moody building $X$, the set 
of fixed points
$\Delta^{\mu(\overline\Lambda_x)}$  is a non-empty closed convex 
union of facets in the Euclidean
building $\Delta$.

{\rm (ii)~} If the image $\mu(\overline\Lambda)$ is non-trivial, then
for any pair of distinct vertices $v \neq v'$ in $X$, the sets of fixed points
$\Delta^{\mu(\overline\Lambda_v)}$ and 
$\Delta^{\mu(\overline\Lambda_{v'})}$ are disjoint.
\rm\pec

{\sc Remark}.---~
As $W$ is infinite, (ii) implies that $\mu(\overline\Lambda)$ is 
unbounded in ${\bf
G}(k)$ since it has no global fixed point in $\Delta$.

\pec

{\it Proof}.---~ 
(i). Let $x \! \in \! X$.
By the very definition of the topology on ${\rm Aut}(X)$, the group 
$\overline\Lambda_x$ is compact, and so is its
continuous image $\mu(\overline\Lambda_x)$.
By non-positive curvature of $\Delta$ and the Bruhat-Tits fixed point 
lemma \cite[VI.4]{Bro89},
the set of $\mu(\overline\Lambda_x)$-fixed points in $\Delta$ is non-empty.
The ${\bf G}(k)$-action on $\Delta$ is simplicial, and since 
$\overline\Lambda$ is
topologically simple (\ref{sss - top simple}) and $\mu$ is continuous, the
$\mu(\overline\Lambda)$-action is type-preserving.
Therefore each time a subgroup of $\mu(\overline\Lambda)$ fixes a 
point, it fixes the closure
of the facet containing it.
The convexity of $\Delta^{\mu(\overline\Lambda_x)}$ comes from the 
intrinsic definition of a
geodesic segment in $\Delta$
\cite[VI.3A]{Bro89}, and from the fact that ${\bf G}(k)$ acts 
isometrically on $\Delta$.

\pec

(ii). The type of a vertex in the ${\rm CAT}(0)$-realization of a 
building defines a subdiagram in the
Dynkin diagram of $\Lambda$ which is spherical and maximal for this 
property \cite[4.3]{Rem99}.
Therefore the fixator of a vertex $v$ in the building $X$ is a 
maximal spherical parabolic subgroup of
the Tits system of \ref{sss - refined TS} with Borel subgroup the 
Iwahori subgroup ${\cal B}$.
Now let $v'$ be another vertex in $X$.
By \cite[IV.2.5]{Bou81}, the group generated by $\overline\Lambda_v$ 
and $\overline\Lambda_{v'}$
is a parabolic subgroup of the latter Tits system, which cannot be 
spherical since it is strictly
bigger than $\overline\Lambda_v$.   By the second remark in \ref{sss 
- Levi factors}, the closed
subgroup
$\overline{\langle \overline\Lambda_v,\overline\Lambda_{v'} \rangle}$ 
generated by $\overline\Lambda_v$ and
$\overline\Lambda_{v'}$ contains an infinite topologically simple 
group $H$; up to conjugacy, this group is a factor of
Proposition \ref{sss - Levi factors} (iv).

\pec

We assume now that there exist two vertices $v \neq v'$ in $X$ such 
that $\Delta^{\mu(\overline\Lambda_v)} \cap 
\Delta^{\mu(\overline\Lambda_{v'})}$
contains a point $y \! \in \! \Delta$.
Then the image of $\overline{\langle 
\overline\Lambda_v,\overline\Lambda_{v'} \rangle}$ by $\mu$ lies in 
the compact fixator of $y$ in ${\bf G}(k)$.
Restricting $\mu$ to the topologically simple group $H$, and 
composing with an embedding
of $k$-algebraic groups ${\bf G} < {\rm GL}_m$, we are in position to 
apply the lemma below: $\mu(H)$ is trivial since it is
bounded.
But then the kernel of $\mu$ is non trivial hence equal to 
$\overline\Lambda$ by topological simplicity (\ref{sss - top simple}).
Finally $\mu(\overline\Lambda) \neq \{ 1 \}$ implies
$\Delta^{\mu(\overline\Lambda_v)} \cap 
\Delta^{\mu(\overline\Lambda_{v'})}=\varnothing$ whenever $v \neq v'$.
\qed\pec

We finally prove the following quite general and probably well-known lemma.

\pec

{\sc Lemma}.---~\it
Let $H$ be an infinite topological group all of whose proper closed 
normal subgroups are finite.
Let $k$ be a non-Archimedean local field and let $\mu: H \to {\rm 
GL}_m(k)$ be a continuous homomorphism for some
$m \geq 1$.
Then $\mu(H)$ is either trivial or unbounded.
\rm\pec

{\it Proof}.---~We denote by ${\cal O}$ the valuation ring and chose 
a uniformizer $\varpi$.
Let us assume that $\mu(H)$ is bounded.
After conjugation, we may -- and shall -- assume that $\mu(H) < {\rm 
GL}_m({\cal O})$
\cite[1.12]{PR94}.  For each integer $N \geq 1$, the group $K:={\rm 
GL}_m({\cal O})$ has an open
finite index congruence subgroup
$K(N) \triangleleft_{\rm f.i.} K$, by definition
$\ker \bigl( {\rm GL}_m({\cal O}) \to {\rm GL}_m({\cal O}/\varpi^N) \bigr)$.
For each $N \geq 1$, $\mu^{-1}\bigl( K(N) \bigr)$ is a closed normal 
finite index
subgroup of $H$.  Since $H$ is infinite, so is $\mu^{-1}\bigl( K(N) 
\bigr)$, and by the
hypothesis on closed normal subgroups of $H$ we have
$H=\mu^{-1}\bigl( K(N) \bigr)$.
Since $\bigcap_{N \geq 1} K(N) = \{ 1 \}$, we have $\mu(H)=\{1\}$.
\qed
\gec

\mysection{Embedding theorem}
\label{s - Embedding theorem}

We use the commensurator super-rigidity to prove that the linearity 
of a countable Kac-Moody group
implies that the corresponding topological group is a closed subgroup 
of a non-Archimedean Lie group.
More precisely:

\pec

{\sc Theorem}.---~\it
Let $\Lambda$ be an almost split Kac-Moody group over the finite 
field ${\bf F}_q$ of
characteristic $p$ with $q > 4$ elements, with infinite Weyl group 
$W$ and buildings $X$ and
$X_-$.
Let $\overline\Lambda$ be the corresponding Kac-Moody topological group.
We make the following assumptions:

\pec

\hskip 4mm {\rm (TS)~} the group $\overline\Lambda$ is topologically simple;

\hskip 4mm {\rm (NA)~} the group $\overline\Lambda$ is not amenable;

\hskip 4mm {\rm (LT)~} the group $\Lambda$ is a lattice of $X \times 
X_-$ for its diagonal
action.

\pec

Then, if $\Lambda$ is linear over a field of characteristic $p$, there exists:

-- a local field $k$ of characteristic $p$ and a connected adjoint $k$-simple
group ${\bf G}$,

-- a topological embedding $\mu: \overline\Lambda \to {\bf G}(k)$ 
with  Hausdorff unbounded
and Zariski dense image,

  -- and a $\mu$-equivariant embedding $\iota: {\cal V}_X \to {\cal V}_\Delta$
from the set of vertices of the Kac-Moody building $X$ of $\Lambda$ 
to the set of vertices of
the Bruhat-Tits building $\Delta$ of ${\bf G}(k)$.
\rm\pec

We now discuss the hypotheses.
When $q>4$ and $\Lambda$ is generated by its root groups, condition 
(TS) is equivalent to
the connectedness of the Dynkin diagram of $\Lambda$ (\ref{sss - top simple}).
Condition (LT) is equivalent to $\Gamma$ being a lattice of
$\overline\Lambda$, which holds whenever $q >\!\! >1$ (\ref{sss - KM 
lattices}).
For condition (NA) we note that $\Lambda$ is chamber-transitive; in 
particular, it is non-compact and
$\overline\Lambda$ is cocompact in ${\rm Aut}(X)$, so that 
$\overline\Lambda$ is amenable (resp. Kazhdan) if
and only if ${\rm Aut}(X)$ is.
Therefore, since amenable groups with property (T) are compact, 
condition (NA) is satisfied whenever the latter group
${\rm Aut}(X)$ has property (T).
By \cite{DJ02} the question of property (T) for automorphism groups 
of buildings admits a
fairly complete answer.
Another case when (NA) is fulfilled is when the building $X$ admits a 
CAT$(-1)$ metric.
Indeed, amenability of ${\rm Aut}(X)$ would imply the existence of a 
probability measure on
$\partial_\infty X$ stabilized by the latter group.
Since ${\rm Aut}(X)$ is non-compact, the relevant Furstenberg lemma 
in this context \cite[Lemma
2.3]{BM96} would imply that the support of this measure consists of 
at most two points.
This would imply the stability of a boundary point or a geodesic, in 
contradiction with the
chamber-transitivity of the $\Lambda$-action on $X$.

\pec

{\sc Example}.---~ 
The above three conditions are satisfied for instance whenever the 
Kac-Moody building $X$ has
apartments isomorphic to a hyperbolic tiling and the thicknesses at 
panels are large enough, a
widely available situation that we will meet in Sect. \ref{s - counter-ex}.

\pec

The remainder of this section is devoted to the proof of the above 
embedding theorem.

\pec

\mysubsection{Semisimple Zariski closure and injectivity.---~} 
\label{ss - Z closures and injectivity}
The linearity assumption says that there is a field ${\bf F}$ of 
characteristic $p$ and an injective group homomorphism
$\eta: \Lambda \hookrightarrow {\rm GL}_N({\bf F})$ to a general 
linear group over ${\bf F}$.
We choose an algebraic closure $\overline {\bf F}$ of ${\bf F}$, and 
denote by ${\bf H}$ the
Zariski closure $\overline{\eta(\Lambda)}^Z$ of the image of $\eta$ 
in ${\rm GL}_N$.
Let ${\cal R}{\bf H}^\circ$ be the radical of the identity component 
${\bf H}^\circ$ of
${\bf H}$ in the Zariski topology.
The group ${\bf H}^\circ/{\cal R}{\bf H}^\circ$ is connected normal and of
finite index in ${\bf H}/{\cal R}{\bf H}^\circ$, hence it is its 
(semisimple) identity
component.
We denote by ${\bf L}$ the quotient of ${\bf H}/{\cal R}{\bf H}^\circ$ by its
finite center ${\bf Z}({\bf H}/{\cal R}{\bf H}^\circ)$, and by $q$ the
natural quotient map ${\bf H} \to {\bf L}$.
Note that ${\bf L}^\circ$ is adjoint semisimple and that the identity 
component of
${\rm ker}(q)$ is solvable.
We consider the composed homomorphism:

\pec

\centerline{$\varphi: \,\, \Lambda \,\, {\buildrel \eta \over \to} 
\,\, {\bf H} \,\,
{\buildrel q \over \to} \,\, {\bf L}$.}

\pec

Let ${\bf L} < {\rm GL}_M$ be an embedding of linear algebraic groups.
Since $\Lambda$ is finitely generated, there is a finitely generated 
field extension ${\bf E}
\subset \overline {\bf F}$ of ${\bf Z}/p{\bf Z}$ such that 
$q(\Lambda)<{\rm GL}_M({\bf E})$.
The group $\varphi(\Lambda)$ is Zariski dense in ${\bf L}$, so the group
${\bf L}$ is defined over
${\bf E}$ \cite[AG 14.4]{Bor90}, and we have $\varphi(\Lambda)<{\bf 
L}({\bf E})$.

\pec

{\sc Lemma}.---~\it
The group homomorphism $\varphi : \Lambda \to {\bf L}({\bf E})$ is injective.
\rm\pec

{\it Proof}.---~ 
Let us assume that the kernel of $\varphi$ is non-trivial, so that 
its closure is a non-trivial closed normal
subgroup of $\overline\Lambda$.  By topological simplicity (\ref{sss 
- top simple}) we have:
$\overline{{\rm ker}(\varphi)}=\overline\Lambda$.
But since $\varphi = q \circ \eta$ and since $\eta$ is injective, we have:
${\rm ker}(\varphi) = \eta^{-1} \bigl( \eta(\Lambda) \cap {\rm ker}(q) \bigr)
\simeq \eta(\Lambda) \cap {\rm ker}(q)$.
This would imply that ${\rm ker}(\varphi)$ is virtually solvable, 
hence amenable for the
discrete topology \cite[4.1.7]{Zim84}.
Then $\overline{{\rm ker}(\varphi)}=\overline\Lambda$ and 
\cite[4.1.13]{Zim84}  would imply that
$\overline\Lambda$ is amenable, which is excluded by the assumption (NA).
\qed
\pec

\mysubsection{Unbounded image and continuous extension.---~} 
\label{ss -Unbounded image and continuous extension}
Our next goal is to check that we are in position to apply the 
commensurator super-rigidity theorem.

\pec

{\sc Proposition}.---~\it
{\rm (i)~} There exist an infinite order element $\gamma$ in the 
lattice $\Gamma$ and a field embedding
$\sigma : {\bf E} \to k$ into a local field of characteristic $p$ 
such that $\sigma\bigl(\varphi(\gamma)\bigr)$
is semismple with an eigenvalue of absolute value $>1$ in the adjoint 
representation of ${\bf L}(k)$.

{\rm (ii)~} There is a connected adjoint $k$-simple group ${\bf G}$ 
and an injective continuous group
homomorphism $\mu : \overline\Lambda \to {\bf G}(k)$.
The map $\mu$ coincides on a finite index subgroup of
$\Lambda$ with the composition of $\varphi$ with the projection onto 
a $k$-simple factor of ${\bf L}^\circ$;
its image is Zariski dense and Hausdorff unbounded.
\rm\pec

{\it Proof}.---~ 
(i). Let us fix a reflection $s$ in a wall $H_s$ containing an edge 
of the standard chamber
$R$, and let us denote by $a_s$ the simple root bordered by $H_s$.
The condition (TS) implies that the Weyl group $W$ of $\Lambda$ (or 
$\overline\Lambda$)
is indecomposable (of non-spherical type).
Therefore by \cite[Proposition~8.1 p.~309]{Hee93} there is a 
reflection $r$ in a wall $H_r$ such
that $rs$ has infinite order, implying that $H_r$ and $H_s$ don't 
meet in the interior of the
Tits cone of $W$ \cite[5.2]{Rem99}.
We call $-b$ the negative root bordered by $H_r$.
If $-a_s \cap -b =\varnothing$ then by definition $\{ -a_s;-b\}$ is a 
non-prenilpotent pair of
roots \cite[3.2]{Tit87}; otherwise $\{ -a_s;-s.b\}$ is.
In any case the group generated by the corresponding root groups is
isomorphic to
${\bf F}_q * {\bf F}_q$ \cite[Proposition 5]{Tit90}.
Therefore it contains an element $\tilde\gamma \! \in \! \Gamma$ of 
infinite order, and by
injectivity $\varphi(\tilde\gamma)$ has infinite order too.
Since the field ${\bf E}$ has characteristic $p>0$, a suitable power 
$p^r$ kills the
unipotent part of the Jordan decomposition of $\varphi(\tilde\gamma)$.
Let us set $\gamma:=\tilde\gamma^{(p^r)}$, so that $\varphi(\gamma)$ 
is semisimple.
It has an eigenvalue $\lambda$ of infinite multiplicative order, so 
by \cite[Lemma 4.1]{Tit72}
there is a local field $k$ endowed with a valuation $v$ and a field
embedding $\sigma: {\bf E}[\lambda] \to k$ such that 
$v\bigl(\sigma(\lambda)\bigr) \neq 0$.
Up to replacing $\gamma$ by $\gamma^{-1}$, this proves (i).

\pec

(ii). By (i) the composed map $({\bf L}(\sigma) \circ \varphi) :
\Lambda \to {\bf L}({\bf E}[\lambda]) \to {\bf L}(k)$, which for 
short we still denote by $\varphi$,
is such that $\Gamma$ has unbounded image in ${\bf L}(k)$.
Let us introduce the preimage $\Lambda^\circ :=\varphi^{-1}({\bf L}^\circ)$.
It is a normal finite index subgroup of $\Lambda$, which we denote by
$\Lambda^\circ \triangleleft_{\rm f.i.} \Lambda$.
We also set $\Gamma^\circ := \Gamma \cap \Lambda^\circ$.
Since $\Gamma^\circ \triangleleft_{\rm f.i.} \Gamma$, 
$\varphi(\Gamma^\circ)$  is not relatively compact in
${\bf L}^\circ(k)$.
The connected adjoint semisimple $k$-group ${\bf L}^\circ$ decomposes 
as a direct product of adjoint
connected $k$-simple factors.
One of them, which we denote by ${\bf G}$, is such that the 
projection of $\varphi(\Gamma^\circ)$ is not
relatively compact.
The abstract group homomorphism we consider now, and which we denote by
$\varphi \! \mid_{\Lambda^\circ}$, is obtained by composing with the 
projection onto ${\bf G}$.
Therefore we obtain $\varphi \! \mid_{\Lambda^\circ}: \Lambda^\circ 
\to {\bf G}(k)$ such that
$\varphi(\Gamma^\circ)$ is unbounded in ${\bf G}(k)$.
We also have: $\overline{\varphi(\Lambda^\circ)}^Z={\bf G}$.

\pec

By Lemma \ref{sss - KM lattices} the group $\Lambda$ is contained in 
the commensurator
${\rm Comm}_{\overline\Lambda}(\Gamma)$.
Since $\Gamma^\circ <_{\rm f.i.} \Gamma$, we have:
${\rm Comm}_{\overline\Lambda}(\Gamma)={\rm 
Comm}_{\overline\Lambda}(\Gamma^\circ)$, so we are in
position to apply the commensurator superrigidity theorem of the 
Appendix in order to extend
$\varphi \! \mid_{\Lambda^\circ}$ to a continuous homomorphism
$\mu: \,\, \overline{\Lambda^\circ} \,\, \to \,\, {\bf G}(k)$,  where 
$\overline{\Lambda^\circ}$ denotes the closure of
$\Lambda^\circ$ in ${\rm Aut}(X)$.
The non-trivial closed subgroup $\overline{\Lambda^\circ}$ is normal 
in $\overline\Lambda$,
hence it is $\overline\Lambda$ by topological simplicity (\ref{sss - 
top simple}).
Therefore there is a map $\mu: \,\, \overline\Lambda \,\, \to \,\, 
{\bf G}(k)$ which coincides with the abstract
group homomorphism $\varphi$ on $\Lambda^\circ$.
By topological simplicity of $\overline\Lambda$, $\mu$ is either 
injective or trivial.
By Zariski density of the image, the only possible case is that $\mu$ 
be injective.
Summing up, we have obtained an injective continuous group homomorphism
$\mu: \,\, \overline\Lambda \,\, \hookrightarrow \,\, {\bf G}(k)$ 
such that $\mu(\Gamma^\circ)$ is unbounded
in ${\bf G}(k)$ and $\overline{\mu(\Lambda^\circ)}^Z={\bf G}$.
\qed\pec

\mysubsection{Embedding of vertices and closed image.---~} 
\label{ss - Embedding of vertices and closed image}
We can finally conclude in view of the following lemma.

\pec

{\sc Lemma}.---~\it
{\rm (i)~} There is a $\mu$-equivariant injective unbounded map
$\iota: {\cal V}_X \hookrightarrow {\cal V}_\Delta$ from the vertices 
of the Kac-Moody building $X$ into the
vertices of the Bruhat-Tits building $\Delta$ of ${\bf G}(k)$.

{\rm (ii)~} The continuous homomorphism $\mu$ sends closed subsets of 
$\overline\Lambda$ to closed subsets
of ${\bf G}(k)$.
\rm\pec

This is the part of the proof which most uses Kac-Moody and Tits 
system theories, so let us briefly recall some facts.
See \cite[\S 5]{Ron89}  for the general connection between Tits 
systems and buildings,
\ref{ss - Tits system} for our specific case.
We keep the inclusion $R \subset A$ of the standard chamber in the 
standard apartment.
Let $W_R$ be the quotient of the stabilizer $N_A={\rm 
Stab}_{\overline\Lambda}(A)$ by the fixator
$\Omega_A={\rm Fix}_{\overline\Lambda}(A)$, which is the Weyl group 
of the building $X$, and of the groups
$\Lambda$ and $\overline \Lambda$.
It is generated by the reflections along the panels of $R$, and 
simply transitive on the chambers of $A$.
We denote by ${\cal B}$ the standard Iwahori subgroup $\overline\Lambda_R$.
By Theorem \ref{sss - refined TS}, ${\cal B}$ is the Borel subgroup 
of a Tits system in $\overline\Lambda$ with Weyl group
$W_R$.
The Tits system structure implies a Bruhat decomposition \cite[IV.2.3]{Bou81}:
$\overline\Lambda=\bigsqcup_{w \in W_R} {\cal B}w{\cal B}$.

\pec

{\it Proof}.---~
(i). By \ref{sss - maps} (i) we choose for each vertex $v$ in the 
closure of the chamber $R$ a
$\mu(\overline\Lambda_v)$-fixed vertex $\iota(v) \! \in \! \Delta$.
We can extend  $\mu$-equivariantly this choice
$\overline R \cap {\cal V}_X \to {\cal V}_\Delta$ to obtain  a map 
$\iota : {\cal V}_X \to
{\cal V}_\Delta$, where $\iota(v)$ is a 
$\mu(\overline\Lambda_v)$-fixed vertex in $\Delta$ for
each vertex $v$ in $X$.  By \ref{sss - maps} (ii), the sets of fixed points
$\Delta^{\mu(\overline\Lambda_v)}$ are mutually disjoint when $v$ 
ranges over ${\cal V}_X$, so
$\iota$ is injective.
By discreteness of the vertices in $\Delta$, the unboundedness of 
$\iota$ follows from its
injectivity because ${\cal V}_X$ is infinite (since so is $W$).

\pec

(ii). Let $F$ be a closed subset of $\overline\Lambda$; we must show 
that $\overline{\mu(F)} < \mu(F)$.
Let $\displaystyle g = \lim_{n \to \infty} \mu(h_n)$ be in 
$\overline{\mu(F)}$, with
$h_n \! \in \! F$ for each $n \geq 1$.
It is enough to show that $\{ h_n \}_{n \geq 1}$ has a converging subsequence.
By the Bruhat decomposition $\overline\Lambda = \bigsqcup_{w \in W_R} 
{\cal B}w{\cal B}$, we can write
$h_n=k_nw_nk_n'$ with $k_n,k_n' \! \in \! {\cal B}$ and $w_n \! \in \! N_A$.
Since by compactness of ${\cal B}$ the sequences $\{ k_n \}_{n \geq 
1}$ and $\{ k_n' \}_{n \geq 1}$ admit
cluster values, we are reduced to the situation where
$\displaystyle g = \lim_{n \to \infty} \mu(w_n)$ with $w_n \! \in \! 
N_A$ for each
$n \geq 1$.

\pec

Let us assume that the union of chambers $\bigcup_{n \geq 1} w_n.R$ 
is unbounded in $A$.
Then there is an injective subsequence of chambers $\{ w_{n_j}.R 
\}_{j \geq 1}$.
Let us fix a vertex $v \! \in \! \overline R$.
Since its stabilizer in $W_R$ is finite, possibly after extracting 
again a subsequence,
we get an injective sequence of vertices $\{ w_{n_j}.v \}_{j \geq 1}$.
But $\mu(w_{n_j}).\iota(v)=\iota(w_{n_j}.v)$ where $\iota: {\cal V}_X 
\hookrightarrow {\cal V}_\Delta$
is the $\mu$-equivariant embedding of vertices of (i).
Since $\displaystyle g = \lim_{n \to \infty} \mu(w_n)$, the 
continuity of the ${\bf G}(k)$-action
on $\Delta$ implies: $\displaystyle \lim_{j \to \infty} 
\mu(w_{n_j}).\iota(v) = g.\iota(v)$.
By discreteness of the vertices in $\Delta$, the sequence $\{ 
\iota(w_{n_j}.v) \}_{j \geq 1}$ hence the sequence
$\{ w_{n_j}.v \}_{j \geq 1}$ is eventually constant: a contradiction.

\pec

We henceforth know that the sequence $\{ w_n.R \}_{n \geq 1}$ is 
bounded, hence takes finitely many values.
So there is a subsequence $\{ w_{n_j} \}_{j \geq 1}$ and $w \! \in \! 
N_A$ such that
$\{ w_{n_j} \}_{j \geq 1}$ is constant equal to $w$ modulo $\Omega_A$.
This proves the lemma, possibly after extracting a converging 
subsequence in the compact group
$\Omega_A<{\cal B}$.
\qed
\gec

\mysection{Some concrete non-linear examples}
\label{s - counter-ex}

We prove that most of countable Kac-Moody groups with right-angled 
Fuchsian buildings  are not
linear over any field.
This requires to settle structure results for boundary point 
stabilizers and generalized unipotent radicals.
Dynamical arguments due to G. Prasad play a crucial role.
The geometry of compactifications of
buildings sheds some light on ideas of the proof.

\pec

\mysubsection{Groups with right-angled Fuchsian buildings.---~}
\label{ss - Fuchsian}
We prove further structure results on topological Kac-Moody groups 
with right-angled
Fuchsian buildings.
Let $R$ be a regular hyperbolic right-angled $r$-gon, so that $r \geq 5$.
To obtain a countable Kac-Moody group with a building covered by 
chambers $\simeq R$, we need to lift the Coxeter
diagram of the corresponding Fuchsian Weyl group $W_R$ to a Dynkin diagram.
The Coxeter diagram of the latter group  is connected and all its 
edges are labelled by $\infty$, so according to the
rule \cite[3.1]{Tit87}  infinitely many Dynkin diagrams are suitable.
Henceforth, $\Lambda$ denotes a Kac-Moody group over ${\bf F}_q$ whose positive
building $X$ is isomorphic to some $I_{r,1+q}$ (\ref{sss - uniform lattices}).
We choose a standard positive chamber $R$ in a standard positive 
apartment $A \simeq {\bf H}^2$
(\ref{sss - countable Tits}).
We denote by $d$ the natural ${\rm CAT}(-1)$ distance on $X$ and by 
$\ell$ the length of any edge.
We fix a numbering $\{ E_i \}_{i \in {\bf Z}/r}$ by ${\bf Z}/r$ of 
the edges of $R$, and we denote by
$a_i$ the simple root containing $R$ whose wall contains $E_i$.
Note that since all wall intersections are orthogonal, all the edges 
in a given wall $L$ have the same
type, which we also call the {\it type~} of $L$ \cite[4.A]{RR02}.

\pec

\mysubsubsection{} 
\label{sss - para horo} 
A {\it geodesic ray~} in a geodesic ${\rm CAT}(-1)$ metric space 
$(X,d)$ is an isometry
$r:[0;\infty) \to X$.
The {\it Busemann function~} of $r$ is the function $f_r:X\to {\bf 
R}$ defined by
$\displaystyle f_r(x):=\lim_{t \to \infty} \bigl( d(x,r(t))-t \bigr)$.
The {\it horosphere~} (resp. {\it horoball}) associated to $r$ is the level set
$H(r):=\{ f_r=0 \}$ (resp. $Hb(r):=\{ f_r\leq 0 \}$).
Let $\xi \! \in \! \partial_\infty X$ be a boundary point of $X 
\simeq I_{r,1+q}$, i.e. an asymptotic class of geodesic ray
\cite[\S 7]{GH90}.

\pec

{\sc Definition}.---~\it
{\rm (i)~} We call {\rm parabolic subgroup attached to $\xi$~} the 
stabilizer $P_\xi:={\rm Stab}_{\overline\Lambda}(\xi)$.

{\rm (ii)~} We call {\rm horospheric subgroup attached to $\xi$~} the 
subgroup of $P_\xi$ stabilizing each
horosphere centered at $\xi$.
We denote it by $D_\xi$.
\rm\pec

{\sc Remarks}.---~
1) The terminology mimicks the geometric definition of proper 
parabolic subgroups in semisimple Lie groups, but
there are differences.
There are two kinds of boundary points, according to whether the 
point is the end
of a wall or not.  A point of the first kind is called {\it 
singular}; otherwise, it is called {\it regular}.
This provides two kinds of parabolic subgroups which are both 
amenable by \cite[Proposition
1.6]{BM96}.  The only classical case when all proper parabolic 
subgroups are amenable is when
they are minimal, i.e. when the split rank equals one.  But then all 
proper parabolics are
conjugate.

2) Defining a suitable notion of rank for the buildings $I_{r,1+q}$ 
is an interesting question.
On the one hand, they contain sharply different kinds of lattices, 
which makes them close to trees
(\ref{sss - uniform lattices}).
On the other hand, they enjoy remarkable rigidity properties, which 
makes them close to higher-rank buildings
\cite{BP00}.

\pec

{\sc Picture}.---~

\pec
\centerline{\begin{picture}(0,0)%
\includegraphics{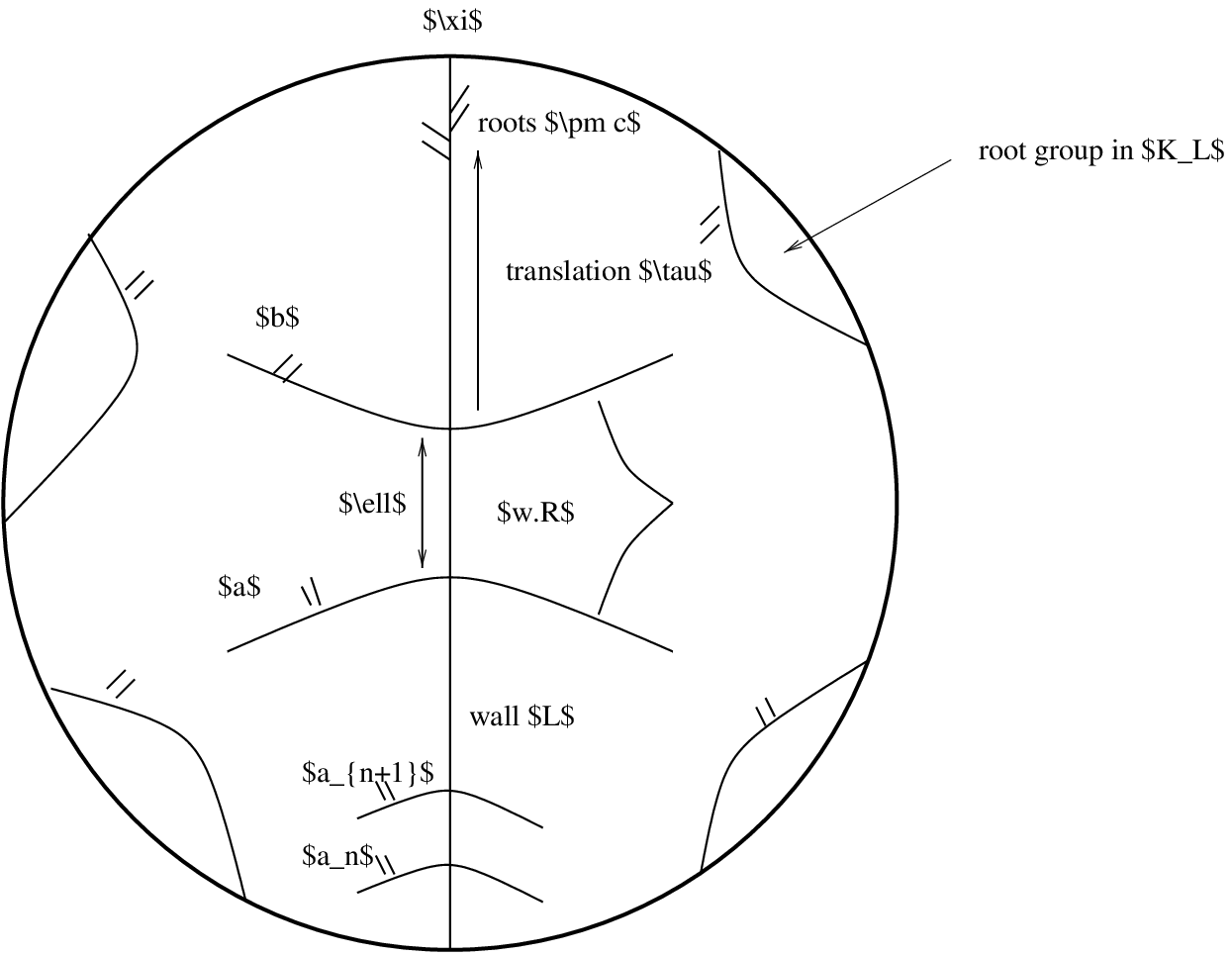}%
\end{picture}%
\setlength{\unitlength}{1973sp}%
\begingroup\makeatletter\ifx\SetFigFont\undefined
\def\x#1#2#3#4#5#6#7\relax{\def\x{#1#2#3#4#5#6}}%
\expandafter\x\fmtname xxxxxx\relax \def\y{splain}%
\ifx\x\y   
\gdef\SetFigFont#1#2#3{%
  \ifnum #1<17\tiny\else \ifnum #1<20\small\else
  \ifnum #1<24\normalsize\else \ifnum #1<29\large\else
  \ifnum #1<34\Large\else \ifnum #1<41\LARGE\else
     \huge\fi\fi\fi\fi\fi\fi
  \csname #3\endcsname}%
\else
\gdef\SetFigFont#1#2#3{\begingroup
  \count@#1\relax \ifnum 25<\count@\count@25\fi
  \def\x{\endgroup\@setsize\SetFigFont{#2pt}}%
  \expandafter\x
    \csname \romannumeral\the\count@ pt\expandafter\endcsname
    \csname @\romannumeral\the\count@ pt\endcsname
  \csname #3\endcsname}%
\fi
\fi\endgroup
\begin{picture}(7910,7699)(2366,-6995)
\put(6451,-1561){\makebox(0,0)[lb]{\smash{\SetFigFont{8}{9.6}{rm}translation $\tau$}}}
\put(10276,-586){\makebox(0,0)[lb]{\smash{\SetFigFont{8}{9.6}{rm}root group in $K_L$}}}
\put(5776,464){\makebox(0,0)[lb]{\smash{\SetFigFont{8}{9.6}{rm}$\xi$}}}
\put(4426,-1936){\makebox(0,0)[lb]{\smash{\SetFigFont{8}{9.6}{rm}$b$}}}
\put(5101,-3436){\makebox(0,0)[lb]{\smash{\SetFigFont{8}{9.6}{rm}$\ell$}}}
\put(6376,-3511){\makebox(0,0)[lb]{\smash{\SetFigFont{8}{9.6}{rm}$w.R$}}}
\put(4126,-4111){\makebox(0,0)[lb]{\smash{\SetFigFont{8}{9.6}{rm}$a$}}}
\put(6151,-5161){\makebox(0,0)[lb]{\smash{\SetFigFont{8}{9.6}{rm}wall $L$}}}
\put(6226,-361){\makebox(0,0)[lb]{\smash{\SetFigFont{8}{9.6}{rm}roots $\pm c$}}}
\put(4801,-5611){\makebox(0,0)[lb]{\smash{\SetFigFont{8}{9.6}{rm}$a_{n+1}$}}}
\put(4801,-6286){\makebox(0,0)[lb]{\smash{\SetFigFont{8}{9.6}{rm}$a_n$}}}
\end{picture}
}

\pec

We are now interested in singular boundary points.
Let $L$ be a wall and let $a$ and $b$ be two roots with $a \supset b$ 
whose walls $\partial a$ and $\partial b$
intersect $L$.
The reflection along $\partial a$ (resp. $\partial b$) is denoted by 
$r_a$ (resp. $r_b$) and
$\tau:=r_b.r_a$ is a hyperbolic translation along $L$ with attracting 
point $\xi$ contained in $a$.  We
assume henceforth that the strip $a \cap (-b)$ of $A$ doesn't contain 
any other wall intersecting $L$.  Then
$\tau$ is a {\it two-step hyperbolic translation}, i.e. a translation 
of (minimal) length $2\ell$.
We set: $a_0:=a$, $a_1:=b$, $a_{2j}:=\tau^j.a$, $a_{2j+1}:=\tau^j.b$ 
for $j \! \in \! {\bf Z}$,
so that $a_k \supset a_{k+1}$, and for each $n \! \in \! {\bf Z}$ we 
denote by $v_n$ the vertex $\partial a_n \cap L$.
The image of the geodesic ray $r : [0;\infty) \to X$ such that 
$r(k\ell)=v_{n+k}$ for each $k \geq 0$, is
$[v_n\xi)=a_n \cap L$.

\pec

{\sc Definition}.---~\it
{\rm (i)~} We denote by $V_n$ the closed  subgroup generated by the 
root groups $U_{a_k}$ for $a_k \supset
a_n$, i.e. $k \leq n$.
We define $V_{\xi,A} := \bigcup_{n \in {\bf Z}} V_n$, which we call a 
{\rm generalized unipotent
radical}.

{\rm (ii)~} We denote by $K_L$ the fixator of the wall $L$ in 
$\overline\Lambda$.

{\rm (iii)~} We denote by $M_{L,A}$ the finite reductive group over 
${\bf F}_q$ generated by the torus $T$ and the two
opposite root groups $U_{\pm c}$ such that $\partial c=L$.
\rm\pec

{\sc Remark}.---~
For each $n\geq 1$ the group $V_n$ is compact since it fixes the 
half-space $a_n$ of $A$.

\pec

Two half-walls define the same {\it germ~} if they intersect along a half-wall.
Two half-walls in the same germ clearly define the same boundary 
point, and the converse is true because
any two disjoint closed edges are at distance $\geq\ell$, so we can 
talk of the germ of a singular boundary point.
Since any element of $\overline\Lambda$ sends a wall onto a wall, we 
have the following characterizations.

\pec

{\sc Lemma}.---~\it
The group $P_\xi$ is the stabilizer of the germ of $\xi$,
and $D_\xi$ is the fixator of this germ, meaning that there is a 
half-wall in it which is fixed under $D_\xi$.
\qed\rm\pec

\mysubsubsection{} 
\label{sss - horoballs} 
Thanks to the Moufang property and the language of horoballs, we can 
say more about the groups $V_n$.
Let us denote by $E$ the intersection of the wall $L$ with the strip 
$a \cap (-b)$.
By the previously assumed minimality of $a \cap (-b)$, $E$ is reduced 
to the edge of a chamber $w.R$.
Transforming the objects above by $w^{-1} \! \in \! W_R$,
we may -- and shall -- assume that we are in the case where $L$ is 
the wall $\partial a_i$
(where $i$ the type of $E$), and either
$a=a_{i-1}$ and $b=-a_{i+1}$, or $a=a_{i+1}$ and $b=-a_{i-1}$.
These two situations are completely analogous,
and we assume that $a=a_{i-1}$ and $b=-a_{i+1}$.
We denote by $r_j$ the reflection in the edge $E_j$ of $\overline R$, we set
$J:=\{ i-1;i+1 \}$, $W_J := \langle r_{i-1},r_{i+1} \rangle$
and we use notions and notation of \ref{sss - Levi factors}.
Then $V_{\xi,A}$ is a subgroup of the topologically simple group 
$\overline G_J$,
and the inessential building $X_J$ is a combinatorial tree.
Its  vertices (resp. edges) are the $\overline G_J$-transforms of a 
line $\partial a_{i-1}$ or
$\partial a_{i+1}$ (resp. of the strip $a_{i-1} \cap a_{i+1}$).
The root groups in $\overline G_J$ are those indexed by the roots 
$a=w.a_j$ with
$w \! \in \! W_J$ and $j \! \in \! J$.
They are automorphisms of the tree $X_J$ fulfilling the Moufang 
condition \cite[\S 6.4]{Ron89}.
Selecting in each strip of $X_J$ the $\overline G_J$-transform of the 
edge $E_i$ of $R$
provides a bijection between the inessential tree $X_J$ and the {\it 
tree-wall~}
attached to $E_i$ \cite[1.2.C]{BP00}.
We henceforth adopt the tree-wall viewpoint when dealing with $X_J$,
so that $v_i$ are vertices  and $\xi$ is a boundary point of it.
Moreover the alluded to below geodesic rays, horospheres and 
horoballs are those defined in the tree-wall $X_J$.

\pec

{\sc Lemma}.---~\it
We assume that all the root groups indexed by the roots $a=w.a_j$ 
with $w \! \in \! W_J$ and
$j \! \in \! J$, and containing $\xi$, commute with one another.

\pec

{\rm (i)~} The group $V_n$ acts trivially on the horoball defined by 
$[v_n\xi)$.

{\rm (ii)~} For each vertex $v$ on the horosphere defined by 
$[v_n\xi)$, the group $U_{a_n}$ is
simply transitive on the edges containing $v$ which are outside the 
corresponding horoball.

{\rm (iii)~} We have: $\bigcap_{n \geq 1} V_n = \{ 1 \}$.

{\rm (iv)~} Any $g \! \in \! V_{\xi,A}$ stabilizing $[v_n\xi)=a_n\cap 
L$ belongs to $V_n$.
\rm\pec

{\sc Remark}.---~
The roots $a=w.a_j$ with $w \! \in \! W_J$ and $j \! \in \! J$ are 
the real roots of a
rank 2 Kac-Moody root system, and the corresponding root groups 
generate a rank 2 countable Kac-Moody group
with generalized Cartan matrix $\pmatrix{2&A_{i-1,i+1}\cr A_{i-1,i+1}&2}$.
As explained in \ref{ss - Fuchsian}, these off-diagonal coefficients 
are $\leq -1$, and their product is
$\geq 4$.
According to the explicit commutator relations due to J. Morita 
\cite[\S 3 (6)]{Mor88}, the group generated
by the root groups $U_a$ for which $\bar a \supset \xi$ is abelian 
whenever the off-diagonal coefficients
are both $\leq -2$ (otherwise it may be metabelian), so the 
assumption made in the lemma is quite not
restrictive.

\pec

{\it Proof}.---~ Let $v$ be a vertex in $X_J$ and let $v_N$ be the 
projection of $v$ on the
geodesic $\{ v_i \}_{i \in {\bf Z}}$ for some $N \! \in \! {\bf Z}$.
By the Moufang property, there are uniquely defined $m \geq 0$ and 
$u_i \! \in \! U_{a_i}$ for
$N-m< i \leq N$ such that $v=(u_N u_{N-1}...u_{N-m+1}).v_{N-m}$.
Denoting by $f_\rho$ the Busemann function of $\rho=[v_0\xi)$, we have:
$f_\rho(v)=(m-N)\ell$.
The Moufang property thus provides a parametrization of the horoballs 
centered at $\xi$ in $X_J$ since
$H([v_j\xi))= \{ (u_N u_{N-1}...u_{N-m+1}).v_{N-m} : N \! \in \! {\bf 
Z}$, $m \geq 0$,
$N-m=j$, $u_N \neq 1$ and $u_i \! \in \! U_{a_i}$ for $N-m< i \leq N \}$.
The commutation of all the root groups $U_{a_k}$, along with the 
parametrization of the horoballs by means
of root groups, proves (i) and (ii).
Moreover if $h \! \in \! \bigcap_{n \geq 1} V_n$,
then (i) implies that
$h$ fixes all the horoballs centered at $\xi$ in $X_J$, so $h$ 
belongs to the kernel of the $\overline
G_J$-action on $X_J$.
We have $h=1$ by topological simplicity (\ref{sss - top simple}), 
which  proves (iii).

\pec

(iv).
By (i) any $g \! \in \! V_{\xi,A}$ stabilizes the horospheres centered
at $\xi$, so if $g$ stabilizes $[v_n\xi)$, it fixes it.
By definition of $V_{\xi,A}$ as an increasing union of groups,
it is enough to show that ${\rm Fix}_{V_N}([v_n\xi))=V_n$ for each $N>n$.
Let $\{ K_j \}_{j \geq 1}$ be an increasing exhaustion of the tree-wall
$X_J$ by finite unions of closed facets, and such that $K_1=[v_n;v_N]$.
For each $j \geq 1$, the subset $C_j :=V_N.K_j$ is $V_N$-stable by
construction and is still a finite union of closed facets since $V_N$ is
compact.
By definition of $V_N$ as a closure, we have:
${\rm Fix}_{V_N}([v_n\xi))
= \limproj_{j\geq1}
{\rm Fix}_{\langle U_{a_i} : i \leq N \rangle}([v_n\xi)) \! \mid_{C_j}$, and
an element in
${\rm Fix}_{\langle U_{a_i} : i \leq N \rangle}([v_n\xi)) \! \mid_{C_j}$
can be written $\prod_{i \leq N} u_i\! \mid_{C_j}$
with $u_i \! \in \! U_{a_i}$ and finitely many non-trivial $u_i$'s.
By (ii), we have $u_i=1$ for $n+1 \leq i \leq N$, which implies:
${\rm Fix}_{V_N}([v_n\xi)) \! \mid_{C_j}
= {\rm Fix}_{\langle U_{a_i} : i \leq n \rangle}([v_n\xi)) \! \mid_{C_j}
= V_n \! \mid_{C_j}$.
Passing to the projective limit shows (iv).
\qed\pec

\mysubsubsection{} 
\label{sss - structure para} 
We can now state the  main properties of the boundary point 
stabilizers and of their generalized unipotent
radicals in the commutative case, keeping the previous notation.

\pec

{\sc Proposition}.---~\it
We assume that all the root groups indexed by the roots $a=w.a_j$ 
with $w \! \in \! W_J$ and
$j \! \in \! J$, and containing $\xi$, commute with one another.

\pec

{\rm (i)~} The group $V_{\xi,A}$ is closed, normalized by $\langle 
\tau \rangle$ but not by $K_L$.
Each group $V_n$ is abelian of exponent $p$, hence so is $V_{\xi,A}$.

{\rm (ii)~} The group $K_L$ is normalized but not centralized by 
$\langle \tau \rangle$.
It admits a semidirect product decomposition $M_{L,A} \ltimes 
\widehat{U}_L$, where $\widehat{U}_L$ is a pro-$p$
group.
In particular, $K_L$ is virtually pro-$p$.

{\rm (iii)~} The following decompositions hold:
$P_\xi=K_L \cdot \langle \tau \rangle \cdot V_{\xi,A}$ and $D_\xi=K_L 
\cdot V_{\xi,A}$, with trivial pairwise intersections
of the factors:
$\langle \tau \rangle \cap K_L = \langle \tau \rangle \cap V_{\xi,A} 
= K_L \cap V_{\xi,A} = \{ 1\}$.
\rm\pec

{\it Proof}.---~
(i). Let $u \! \in \! \overline V_{\xi,A}$.
We write: $\displaystyle u = \lim_{j \to \infty} u_j$ with $u_j \! 
\in \! V_{\xi,A}$ for each $j \geq 1$,
and we have $u \! \in \! P_\xi$.
By Lemma \ref{sss - para horo}, there is an $n \! \in \! {\bf Z}$ 
such that $u$ sends the geodesic ray
$[v_n\xi)$ to  a geodesic ray contained in $L$ and ending at $\xi$.
For $j>\!\!>1$ we have: $u_j([v_n;v_{n+1}])=u([v_n;v_{n+1}])$.
Since $u([v_n;v_{n+1}])$ is an edge in $L$ and since $u_j$ stabilizes 
the horospheres centered at $\xi$, we
have $u_j([v_n;v_{n+1}])=u([v_n;v_{n+1}])=[v_n;v_{n+1}]$ for $j>\!\!>1$.
Therefore $u_j$ fixes $[v_n\xi)$ for $j>\!\!>1$, implying by Lemma 
\ref{sss - horoballs} (iv)
that for $j>\!\!>1$ the elements $u_j$ lie in the compact group $V_n$.
This implies $u \! \in \! V_n < V_{\xi,A}$, hence (i).

\pec

We turn now to the group-theoretic properties of $V_{\xi,A}$.
By assumption all the groups $U_{a_k}$ commute, so the
continuous commutator map $[.,.]$ is trivial on a topologically 
generating set for each $V_n$.
This proves the commutativity of each group $V_n$.
Any of the commuting root groups $U_{a_k}$ is isomorphic to $({\bf 
F}_q,+)$ so replacing
$[.,.]$ by
${.}^p$ shows that each $V_n$ is of exponent $p$.
By definition, $V_{\xi,A}$ is normalized but not centralized by $\tau$.
Pick a root $a$ containing $L$ so that $U_a < K_L$ and
choose $n \! \in \! {\bf Z}$ such that $-a \cap -a_n=\varnothing$, 
i.e. $\{ a;a_n \}$ is not
prenilpotent.
By \cite[Proposition 5]{Tit90}, the free product $U_a * U_{a_n}$ injects in
$\Lambda$, so for any $u \! \in \! U_a \setminus \{ 1 \}$ and $u' \! 
\in \! U_{a_n} \setminus
\{ 1 \}$ the order of $[u,u']$ is infinite whereas it would divide 
$p$ if $K_L$ normalized
$V_{\xi,A}$.

\pec

(ii). If $\Pi$ is a panel in the wall $L$, we have $K_L < 
\overline\Lambda_\Pi=M_\Pi \ltimes
\widehat{U}_\Pi$  (\ref{sss - refined TS}) and  $M_\Pi=M_{L,A}$ by 
the precise version of the Levi
decomposition \cite[Theorem 6.2.2]{Rem99}.
The group $M_{L,A}$ fixes $L$, so $M_{L,A} < K_L$.
Moreover the kernel $\widehat{U}_L:=\widehat{U}_\Pi \cap K_L$ of the 
restricted map $K_L \to M_{L,A}$ is pro-$p$, and we have:
$M_{L,A} \cap \widehat{U}_L=\{1\}$ by the same argument as for 
\cite[Lemma 1.C.5 (i)]{RR02}.
The group $K_L={\rm Fix}_{\overline\Lambda}(L)$ is normalized by 
$\langle \tau \rangle$ because $\tau$ stabilizes $L$.
Now we pick a root $a \neq \pm c$ with $a \supset L$, so that $U_a<K_L$.
For $M > \!\! > 1$, the root $\tau^M.a$ contains $L$ but $a \cap 
\tau^M.a$ is a strip in $A$, implying that $\{
a;\tau^M.a \}$ is not prenilpotent.
As for (i)  the free product $U_a * U_{\tau^M.a}$ injects in
$\Lambda$, so for $u \! \in \! U_a \setminus \{ 1 \}$ the commutator 
$[\tau^M,u]$ has infinite order whereas it
would be trivial if $K_L$ were centralized by $\tau$.

\pec

(iii). Let $r : [0;\infty) \to X$ be the geodesic ray such that 
$r(n\ell)=v_n$ for each $n \geq
0$, and let $g \! \in \! P_\xi$.
By Lemma \ref{sss - para horo}, there are integers $N \geq 1$ and $t 
\! \in \! {\bf Z}$ such that
$(g.r)(n\ell)=r((n+t)\ell)$ for $n \geq N$.
Since the $\overline\Lambda$-action on $X$ is type-preserving, $t$ is 
an even number, say $2m$, and we have:
$(\tau^{-m}g.r)(n\ell)=r(n\ell)$ for each $n \geq N$.
Thus $d:=\tau^{-m}g$ fixes the geodesic ray $[v_N\xi)$, hence belongs 
to $D_\xi$, and we are  reduced to
decompose $D_\xi$.

\pec

Let $d \! \in \! D_\xi$, which fixes a geodesic ray $[v_N\xi)$ by
Lemma \ref{sss - para horo}.
The link of the vertex $v_N$ is complete bipartite, so there is a chamber $R'$
whose closure contains both $d.[v_N;v_{N-1}]$ and an edge $E'$ contained in
the wall $\partial a_N$.
By the Moufang property, there exists $u_N \! \in \! U_{a_N}$ such that
$(u_N^{-1}d).R'$ is the chamber in $A$ whose closure contains both
$[v_N;v_{N-1}]$ and $E'$;
in particular, $u_N^{-1}d$ fixes the geodesic ray $[v_{N-1}\xi)$.
By a downwards induction, for each $m<N$ we pick $u_m \! \in \! 
U_{a_m}$ such that
$u_m^{-1}u_{m+1}^{-1}...u_N^{-1}d$ fixes the geodesic ray $[v_{m-1}\xi)$.
By compactness of $V_N$, the sequence $\{ u_N...u_{m+1}u_m \}_{m<N}$ 
has a cluster value $u \! \in \! V_{\xi,A}$ such that
$u^{-1}d$ fixes the geodesic $L$, hence belongs to $K_L$.
Taking inverses, we proved: $D_\xi=K_L \cdot V_{\xi,A}$, and along 
with the previous paragraph:
$P_\xi=K_L \cdot \langle \tau \rangle \cdot V_{\xi,A}$ since $\tau$ 
normalizes $V_{\xi,A}$ and $K_L$.

\pec

The trivial intersection $\langle \tau \rangle \cap V_{\xi,A} = \{ 1 
\}$ follows from the fact that
$\langle \tau \rangle \simeq {\bf Z}$ whereas any non-trivial element 
in $V_{\xi,A}$ has order $p$,
and $\langle \tau \rangle \cap K_L = \{ 1 \}$ follows from the fact 
that no non-trivial power $\tau^m$ fixes $L$.
An element in $K_L \cap V_{\xi,A}$ lies in any ${\rm 
Fix}_{V_{\xi,A}}([v_n\xi))$ ($n \! \in \! {\bf Z}$),
hence in any $V_n$ by Lemma \ref{sss - horoballs} (iv), so $K_L \cap 
V_{\xi,A} = \{ 1\}$ follows from
(iii) in the same lemma.
\qed\pec

{\sc Remark}.---~
1) Horoball arguments as in \ref{sss - horoballs} show that each 
group $V_n$ is isomorphic to $({\bf F}_q[[t]],+)$ and
that there is an isomorphism $V_{\xi,A} \simeq ({\bf F}_q(\!( t 
)\!),+)$ under which conjugation by $\tau$ corresponds
to multiplication by $t^{-2}$ and the $t$-valuation corresponds to 
the index ${}_n$.

2) Let us denote by $-\xi$ the other end of $L$, so that $L=(-\xi\xi)$.
By definition, $D_\xi \cap D_{-\xi}$ stabilizes $L$ and actually 
fixes it since $D_\xi$ stabilizes the horospheres
centered at $\xi$.
Therefore we have: $D_\xi \cap D_{-\xi}=K_L$.

\pec

\mysubsection{Dynamics and parabolics.---~}
\label{ss - dynamics}
Let us have a dynamical viewpoint on the above groups.
The prototype for parabolics, used in \ref{sss - RA not linear}, is 
G. Prasad's work in the
algebraic group case \cite{Pra77}.

\pec

\mysubsubsection{} 
\label{sss - F boundary}
A first consequence of the existence of many hyperbolic translations 
is the connection with
Furstenberg boundaries -- see \cite[VI.1.5]{Mar90} for a definition, 
where this notion is
simply called a {\it boundary}.
We denote by ${\cal M}^1(\partial_\infty X)$ the space of probability 
measures on
$\partial_\infty X$;  it is compact and metrizable for the weak-$*$ topology.
This subsection owes its existence to discussions with M. Bourdon and 
Y. Guivarc'h.

\pec

{\sc Lemma}.---~\it
The asymptotic boundary $\partial_\infty X$ is a Furstenberg boundary 
for $\overline\Lambda$.
\rm\pec

{\it Proof}.---~
Let us prove that the $\overline\Lambda$-space $\partial_\infty X$ is 
both minimal and strongly proximal
\cite[VI.1]{Mar90}.

\pec

Strong proximality.
Let $\mu \! \in \! {\cal M}^1(\partial_\infty X)$.
Since some unions of walls in $X$ are trees, the set of singular 
points is uncountable.
Therefore there is a hyperbolic translation $\tau$ along a wall whose 
repelling point is not one of the at most
countably many atoms for $\mu$.
By dominated convergence, the sequence converges in ${\cal 
M}^1(\partial_\infty X)$ to the Dirac mass centered at
the attracting point of $\tau$.

\pec

Minimality.
Let $\xi \! \in \! \partial_\infty X$.
We write it $\xi=r(\infty)$ for a geodesic ray $r : [0;\infty) \to X$ 
with $r(0) \! \in \! R$.
For each $n \geq 1$, $r(n)$ is in the closure of a chamber $g_n.R$ 
with $g_n \! \in \!
\overline\Lambda$.  By the Bruhat decomposition 
$\overline\Lambda=\bigsqcup_{w \in W_R} {\cal
B}w{\cal B}$, we have:
$r(n) \! \in \! k_nw_n.R$, hence $k_n^{-1}.r(n) \! \in \! A$.
Let us denote by $r_n$ the geodesic ray in $A\simeq {\bf H}^2$ 
starting at $r(0)$ and passing through
$k_n^{-1}.r(n)$.
By compactness of $\partial_\infty {\bf H}^2 \simeq S^1$ and ${\cal 
B}$, there is an extraction
$\{ n_j \}_{j \geq 1}$ such that $r_{n_j}(\infty)$ converges to some 
$\eta \! \in \! \partial_\infty A$ and
$k_{n_j}$ converges to some $k \! \in \! {\cal B}$ as $j \to \infty$.
Thus in the $\overline\Lambda$-compactification $X \sqcup 
\partial_\infty X$, we have:
$\displaystyle \xi = \lim_{n \to \infty} r(n) = \lim_{j \to \infty} r(n_j)
= \lim_{j \to \infty} k_{n_j}.(k_{n_j}^{-1}.r(n_j)) = k.\eta$,
which proves that $\partial_\infty A$ is a complete set of
representatives for the ${\cal B}$-action on $\partial_\infty X$.
Since the action of the Weyl group $W_R$, a lattice of ${\rm 
PSL}_2({\bf R})$, is minimal on $\partial_\infty A$, we
proved the minimality of the $\overline\Lambda$-action on $\partial_\infty X$.
\qed\pec

{\sc Remark}.--- Note that the group $\overline\Lambda$ admits a 
Furstenberg boundary on which it doesn't act
transitively, whereas any such boundary for a semisimple algebraic 
group is an equivariant image of the maximal flag
variety \cite[\S 5]{BM96}.

\pec

\mysubsubsection{} 
\label{sss - limits}
Iteration of hyperbolic translations along walls also leads to 
computations of limits of later use for the non-linearity
theorem.

\pec

{\sc Proposition}.---~\it
Let $\tau$ be a hyperbolic translation along a wall $L$, with 
attracting point $\xi$.
Let $v$ be a vertex on the wall $L$.

\pec

{\rm (i)~} We have:
$\displaystyle \lim_{n \to \infty} \tau^n \overline\Lambda_v \tau^{-n} = D_\xi$
in the compact metrizable space ${\cal S}_{\overline\Lambda}$ of 
closed subgroups of
$\overline\Lambda$, endowed with the Chabauty topology.

{\rm (ii)~} For any $u \! \in \! V_{\xi,A}$, we have: $\displaystyle 
\lim_{n \to \infty} \tau^{-n} u \tau^n = 1$

{\rm (iii)~} For any $g \! \in \! D_\xi$, the sequence $\{ \tau^{-n} 
g \tau^n \}_{n \geq 1}$ is bounded in
$\overline\Lambda$.
\rm\pec

{\sc Remark}.---~ 
Point (i) about the Chabauty topology on closed subgroups is used in
the final discussion \ref{sss - compactifications}.
Recall that the Chabauty topology on the closed subsets of a 
topological space $S$ is the topology
defined by a sub-base consisting of the subsets $O(K):=\{A$ closed 
$:A\cap K=\varnothing\}$ for all
compact subsets $K$, and $O'(U):=\{A$ closed $:A\cap 
U\neq\varnothing\}$ for all open subsets $U$
\cite[3.1.1]{CEG87}.
This topology is always compact, and when $S$ is Hausdorff, locally 
compact and secound countable, it is
separable and metrizable \cite[3.1.2]{CEG87}.
When $S$ is locally compact, a sequence $\{ A_n \}_{n \geq 1}$ of 
closed subsets converges in the Chabauty
topology to a closed subset $A$ if and only if:

1) Any limit $\displaystyle x=\lim_{k \to \infty} x_{n(k)}$ for an increasing
$\{ n(k) \}_{k \geq 1}$ with $x_{n(k)} \! \in \! A_{n(k)}$ satisfies 
$x \! \in \! A$.

2) Any $x \! \in \! A$ is the limit of a sequence $\{ x_n \}_{n \geq 
1}$ with $x_n \! \in \! A_n$ for each
$n \geq 1$.

This characterization of convergence is referred to as {\it geometric 
convergence~} \cite[3.1.3]{CEG87}; it
implies that for a locally compact group $G$, the subset ${\cal S}_G$ 
of closed subgroups is closed, hence
compact, for the Chabauty topology.

\pec

{\it Proof}.---~
(i). By compactness of ${\cal S}_{\overline\Lambda}$ it is enough to 
show that any
cluster value $D$ of $\{ \tau^n \overline\Lambda_v \tau^{-n} \}_{n 
\geq 1}$ is equal to  $D_\xi$.
In one direction, the very definition of the Chabauty topology 
implies that $D$ contains $K_L$ and $V_{\xi,A}$, hence $D_\xi$ by
Proposition \ref{sss - structure para} (iii).
Indeed, the group $K_L$ lies in $D$ since it fixes all the vertices 
in $L$, hence lies in all the conjugates
$\tau^n \overline\Lambda_v \tau^{-n}$.
The limit group $D$ also contains $V_{\xi,A}$ since for each $m \! 
\in \! {\bf Z}$ there is $N \! \in \! {\bf N}$ such that
$\tau^n.v \! \in \! a_m$, hence $V_m < \tau^n \overline\Lambda_v 
\tau^{-n}$ for any $n \geq N$.

\pec

We are thus reduced to proving that any cluster value $D$ lies in $D_\xi$.
Let $\nu_v \! \in \! {\cal M}^1(\partial_\infty X)$ be a 
$\overline\Lambda_v$-invariant measure such that
the repelling point of $\tau$ is not an atom for it.
By dominated convergence, we have:
$\displaystyle \lim_{n \to \infty} \tau^n{}_* \nu_v = \delta_\xi$, 
where $\delta_\xi$ is the Dirac mass
at $\xi$, so by \cite[Lemma 3]{Gui01}: $D < P_\xi = {\rm 
Stab}_{\overline\Lambda}(\delta_\xi)$.

\pec

Let $g \! \in \! D$, which by the previous paragraph and Proposition 
\ref{sss - structure para} (iii) we write $g = u \tau^N k$ with
$u \! \in \! V_{\xi,A}$, $N \! \in \! {\bf Z}$ and $k \! \in \! K_L$.
We choose this order to forget the factor $k$ when $g$ acts on $v$.
Since $D$ is a limit group, we also have:
$\displaystyle g = \lim_{j \to \infty} \tau^{n_j} k_j \tau^{-n_j}$ 
for a sequence
$\{ k_j \}_{j \geq 1} \subset \overline\Lambda_v$ and integers $n_j 
\to \infty$ as $j \to
\infty$.  Therefore there is an index $J \geq 1$ for which  $j \geq 
J$ implies $(u \tau^N).v=
(\tau^{n_j} k_j \tau^{-n_j}).v$.  Since $u$ stabilizes all the 
horospheres centered at $\xi$,
there is a vertex $z \! \in \! L$ with
$(u\tau^N).v$ and $\tau^N.v$ at same distance from $z$.
We choose $j > \! \! > 1$ to have $d\bigl( \tau^{n_j}.v, (u \tau^N).v \bigr) =
(n_j-N)\delta$, where $\delta$ is the translation length of $\tau$.
But the group $\tau^{n_j} \overline\Lambda_v \tau^{-n_j}$ stabilizes 
the spheres centered at
$\tau^{n_j}.v$, so that $d\bigl( \tau^{n_j}.v, (\tau^{n_j} k_j 
\tau^{-n_j}).v \bigr) =
n_j\delta$.   In order to have $(u \tau^N).v = (\tau^{n_j} k_j 
\tau^{-n_j}).v$, we must have
$N=0$, i.e. $g=uk$: this shows that $D$ lies in $D_\xi$.

\pec

(ii). Let $u \! \in \! V_{\xi,A}$.
Then $u \! \in \! V_m$ for some $m \! \in \! {\bf Z}$.
For each $N \geq 1$, the sequence $\{\tau^{-n} u \tau^n \}_{n \geq 
N}$ lies in the compact group $V_{m-N}$,
so that any cluster value of $\{\tau^{-n} u \tau^n \}_{n \geq 1}$ 
belongs to $V_{m-N}$.
By Lemma \ref{sss - horoballs} (iii), this shows that the only 
cluster value of the sequence
$\{\tau^{-n} u \tau^n \}_{n \geq 1}$ in the compact subset $V_m$ is 
the identity element, which proves (ii).

\pec

(iii). Let $g \! \in \! D_\xi$, which we write $g=ku$ with $k \! \in 
\! K_L$ and $u \! \in \! V_{\xi,A}$ by Proposition
\ref{sss - structure para} (iii).
Since $u \! \in \! V_m$ for some $m \! \in \! {\bf Z}$, we have $\tau^{-n} 
u \tau^n \! \in \! V_m$ for each $n
\geq 1$.
By Proposition \ref{sss - structure para} (i) $K_L$ is normalized by 
$\tau$, so we finally have
$\tau^{-n} g \tau^n  \!\in\! K_L\cdot V_m$ for each $n \geq 1$.
\qed\pec

\mysubsection{Non-linearity in equal characteristic.---~}
\label{ss - non-linear} 
We finally state and prove the non-linearity theorem for some 
countable Kac-Moody groups with hyperbolic buildings.
It applies to an infinite family of groups, the Weyl group of which 
being of arbitrarily large rank.

\pec

\mysubsubsection{} 
\label{sss - RA not linear}
The previous dynamical results from \ref{ss - dynamics}, as well as 
the embedding theorem from Sect. \ref{s -
Embedding theorem},  provide the  main arguments to prove the result below.

\pec

{\sc Theorem}.---~\it
Let $\Lambda$ be a countable Kac-Moody group over a finite field 
${\bf F}_q$ of characteristic $p$, and let
$r$ be an integer $\geq 5$.
We assume that the geometry of $\Lambda$ is a twinned pair of 
right-angled Fuchsian buildings
$I_{r,q+1}$ with $q \geq \max \{r\!-\!2;5\}$, and that a generalized 
unipotent radical of $\overline\Lambda$
is abelian.
Then $\Lambda$ is not linear over any field.
\rm\pec

{\sc Remark}.---~
According to \ref{sss - horoballs}, the assumption of commutativity 
of a generalized unipotent
radical is mild, since it amounts to requiring that for some $i \! 
\in \! {\bf Z}/r$, both
negative coefficients $A_{i-1,i+1}$ and $A_{i+1,i-1}$ be $\leq -2$ 
(their product must always
be $\geq 4$ to have $X \simeq I_{r,q+1}$ -- see \ref{ss - Fuchsian}).

\pec

{\it Proof}.---~
By \cite[Proposition 4.3]{Rem01}, it is enough to disprove linearity 
in equal characteristic.
Let us assume that there is an abstract injective homomorphism from 
$\Lambda$ to a linear
group in characteristic $p$, in order to obtain  a contradiction.
Up to replacing $\Lambda$ by a finite index subgroup, we may -- and 
shall -- assume that
$\Lambda$ is generated by its root groups.

\pec

We first check that we can apply the embedding theorem of Sect. 
\ref{s - Embedding theorem}.
By the comments in the introduction of this Section, (NA) holds
because $X$ is ${\rm CAT}(-1)$ and $\overline\Lambda$ is chamber-transitive.
By \cite[1.C.1 last remark]{RR02} the growth series of the Weyl group is
$\displaystyle W(t) = {(1+t)^2 \over (1 - (r-2)t + t^2)} \! \in \! 
{\bf Z}[[t]]$, and
$W({1 \over q})$ is finite if and only if $q \geq r-2$; so
  (LT) holds by \ref{sss - KM lattices}.
Theorem  \ref{sss - top simple} implies that (TS) is satisfied 
because the Dynkin diagrams leading to the
buildings $I_{r,q+1}$ are connected \cite[13.3.2]{Rem99} and we have 
assumed that $q>4$.
Consequently, we have a closed embedding $\mu : \overline\Lambda \to 
{\bf G}(k)$ of topological groups.

\pec

Let $L$ be a wall with end $\xi$ such that $V_{\xi,A}$ is abelian.
We pick a hyperbolic translation $\tau$ along $L$ with attracting 
point $\xi$ as in \ref{sss - para horo},
and $u \! \in \! V_{\xi,A} \setminus \{ 1 \}$.
We set $B:=({\rm Ad}\circ\mu)(\tau)$ and $Y:=({\rm Ad}\circ\mu)(u)$,
where ${\rm Ad}$ is the adjoint representation of ${\bf G}$.
Then Lemma \ref{sss - limits} (ii) implies that $\{ B^{-i}YB^i \}_{i 
\geq 1}$ contains the identity element in its
closure, so \cite[Lemma II.1.4]{Mar90} says that the element $B$ has 
two eigenvalues with
different absolute values.  We can thus use \cite[Lemma 2.4]{Pra77}, 
which provides us parabolic
subgroups:  by [loc. cit. (i)] and \ref{sss - limits} (iii), there is 
a proper parabolic
$k$-subgroup ${\bf P}_{\mu(\tau)}$ whose
$k$-points $P_{\mu(\tau)}$ contain $\mu(D_\xi)$.
Replacing $\tau$ by $\tau^{-1}$, the attracting point becomes the 
boundary point $-\xi$ such that $(-\xi\xi)=L$.
By [loc. cit. (ii)], the corresponding parabolic subgroup 
$P_{\mu(\tau^{-1})}$ is opposite
$P_{\mu(\tau)}$, so that
$\mu(D_\xi) \cap \mu(D_{-\xi})$ lies in the Levi factor 
$M_{\mu(\tau)}:=P_{\mu(\tau)} \cap P_{\mu(\tau^{-1})}$.
By the second remark following Proposition \ref{sss - structure 
para}, we have $K_L=D_\xi \cap D_{-\xi}$, sowe
finally obtain: $\mu(K_L)<M_{\mu(\tau)}$.

\pec

The contradiction comes when we look at the image $\mu(V_{\xi,A})$.
By [loc. cit. (i)]  and \ref{sss - limits} (ii), it lies in the 
unipotent radical
${\cal R}_uP_{\mu(\tau)}$.
The decomposition $D_\xi=K_L \cdot V_{\xi,A}$ of Proposition \ref{sss 
- structure para} (iii) then implies:
$\mu(K_L)=\mu(D_\xi) \cap M_{\mu(\tau)}$ and $\mu(V_{\xi,A}) = 
\mu(D_\xi) \cap {\cal R}_uP_{\mu(\tau)}$.
But according to Proposition \ref{sss - structure para} (i), the 
group $V_{\xi,A}$ is not normalized by $K_L$.
\qed\pec

\mysubsubsection{} 
\label{sss - compactifications}
Let us give a geometric flavour to the above proof by using the 
framework of group-theoretic compactifications of
buildings \cite{AGR02}.
We keep the notation of the previous proof, choose a $k$-embedding 
${\bf G} < {\rm GL}_r$ of algebraic groups and
still call $\mu$ the composed closed embedding $\mu : 
\overline\Lambda \to {\rm GL}_r(k)$.
Replacing $\tau$ by $\tau^{(p^r)}$ for $r > \! \! > 1$ and taking a 
finite extension which we still denote by $k$, we
may -- and shall -- assume that $t:=\mu(\tau)$ is diagonal with 
respect to a basis $\{ e_i \}_{1 \leq i \leq r}$ of
$k^r$.
We write: $t.e_i=u_i\varpi^{\nu_i(t)}e_i$ where $\varpi$ is the 
uniformizer of $k$, $u_i \! \in
\! {\cal O}^\times$  and
$\nu_i(t) \! \in \! {\bf Z}$.
Composing $\mu$ with a permutation matrix enables to assume that
$\nu_1(t) \leq \nu_2(t) \leq ... \leq \nu_r(t)$.
The basis $\{ e_i \}_{1 \leq i \leq r}$ defines a maximal flat $F 
\simeq {\bf R}^{r-1}$
in the Bruhat-Tits building
$\Delta$ of ${\rm GL}_r(k)$, whose vertices are the homothety classes 
of ${\cal O}$-lattices
$[\bigoplus_i\varpi^{\nu_i}{\cal O}e_i]$ when 
$\underline\nu:=\{\nu_i\}_{1 \leq i \leq r}$ ranges over ${\bf Z}^r$.
We denote by $o$ the origin $[\bigoplus_i{\cal O}e_i]$.
The same use of \cite[Lemma II.1.4]{Mar90} as in \ref{sss - RA not 
linear} shows that there is
$i \! \in \! \{1;2;... r-1\}$ such that $\nu_i(t) < \nu_{i+1}(t)$.
Geometrically, this means that $\{ t^n.o \}_{n \geq 1}$ is a sequence 
of vertices in the Weyl chamber
$\{ \nu_1 \leq \nu_2 ... \leq \nu_r \}$ which goes to infinity, 
staying in the intersection of the
fundamental walls indexed by the indices $i$ for which $\nu_i(t) = 
\nu_{i+1}(t)$.

\pec

Sequences of points staying at given distance from some walls in a 
Weyl chamber while leaving the others typically
converge in Furstenberg compactifications of symmetric spaces \cite{Gui01}.
This is a hint to consider compactifications of Bruhat-Tits buildings 
in our context \cite{AGR02}.
We denote by ${\cal S}_G$ the space of closed subgroups of a locally 
compact metrizable group $G$,
and we endow ${\cal S}_G$ with the compact metrizable Chabauty 
topology \cite[3.1.1]{CEG87}.
The vertices ${\cal V}_\Delta$ in $\Delta$ are seen as the maximal 
compact subgroups in  ${\rm SL}_r(k)$.
Therefore we have ${\cal V}_\Delta \subset {\cal S}_{{\rm SL}_r(k)}$ 
and we can sum up some
results from \cite{AGR02}:

\pec

{\sc Theorem}.---~\it
The above procedure leads to a ${\rm GL}_r(k)$-compactification of 
$\Delta$  where the boundary points are the
following closed subgroups of ${\rm SL}_r(k)$.
Start with a Levi decomposition $M \ltimes U$ of some proper 
parabolic subgroup and select $K<M$ a maximal
compact subgroup.
Then $K \ltimes U$ is a limit group, and any limit group is of this form.
\qed\rm\pec

{\sc Remark}.---~
In higher rank (i.e. for $r \geq 3$), this compactification is not 
the one obtained by
asymptotic classes of geodesic rays.

\pec

Now Proposition \ref{sss - limits} (i) says that
$\displaystyle \lim_{n \to \infty} \tau^n\overline\Lambda_v\tau^{-n} 
= D_\xi$ in ${\cal S}_{\overline\Lambda}$.
It follows from the geometric characterization of convergence in the 
Chabauty topology (\ref{sss -
limits}), that $\mu$ induces an embedding
$\mu : {\cal S}_{\overline\Lambda} \hookrightarrow {\cal S}_{{\rm SL}_r(k)}$.
Applying $\mu$ to the above limit and using the theorem imply:
$\mu(D_\xi)<K\ltimes U$ and $\mu(D_{-\xi})<K\ltimes U^-$, with $U_-$ 
opposite $U$.
This leads to $\mu(K_L)<K$, an important step in the previous proof.

\pec

The comparison of hyperbolic and Euclidean apartments emphasizes a sharp
difference between Fuchsian and affine root systems.
In a Euclidean apartment, there is a finite number of parallelism 
classes of walls, whereas in the
hyperbolic tiling there are arbitrarily large families of roots 
pairwise intersecting along strips.
This explains why there are so many non prenilpotent pairs of roots 
(hence free products
${\bf F}_q*{\bf F}_q$) in the latter case.
This is used to prove that $K_L$ doesn't normalize $V_{\xi,A}$ 
(\ref{sss - structure para}), a key fact for
non-linearity.
Another crucial difference is the dynamics of the Weyl groups on the
boundaries of apartments: in the hyperbolic case, there are 
infinitely many hyperbolic translations with
strong dynamics \cite[\S 8]{GH90}, whereas the finite index 
translation subgroup of a Euclidean
Weyl group acts trivially on the boundary of a maximal flat.
This makes the computation of limit groups easier in the former case 
(\ref{sss - limits}), but
the boundary of the Furstenberg compactification of a Bruhat-Tits 
building has a much richer
group-theoretic structure since it contains compactifications of 
smaller Euclidean buildings
\cite[\S 14]{Lan96}.

\pec

{\sc Picture}.---~ [on the right hand-side: a compactified apartment 
for ${\rm SL}_3$]

\pec
\centerline{\begin{picture}(0,0)%
\includegraphics{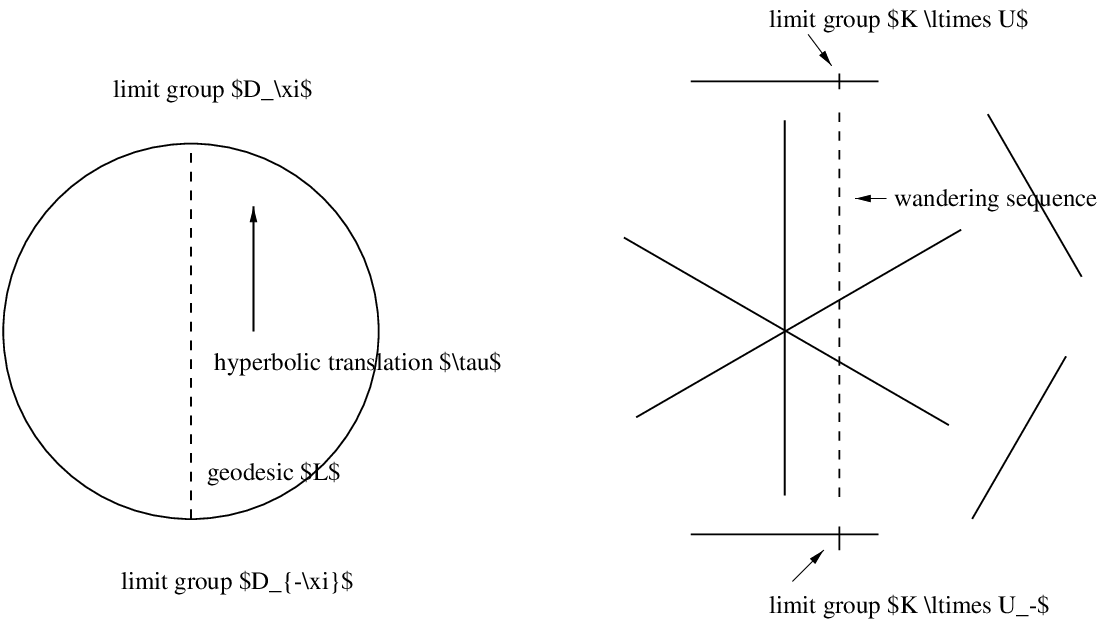}%
\end{picture}%
\setlength{\unitlength}{1973sp}%
\begingroup\makeatletter\ifx\SetFigFont\undefined
\def\x#1#2#3#4#5#6#7\relax{\def\x{#1#2#3#4#5#6}}%
\expandafter\x\fmtname xxxxxx\relax \def\y{splain}%
\ifx\x\y   
\gdef\SetFigFont#1#2#3{%
  \ifnum #1<17\tiny\else \ifnum #1<20\small\else
  \ifnum #1<24\normalsize\else \ifnum #1<29\large\else
  \ifnum #1<34\Large\else \ifnum #1<41\LARGE\else
     \huge\fi\fi\fi\fi\fi\fi
  \csname #3\endcsname}%
\else
\gdef\SetFigFont#1#2#3{\begingroup
  \count@#1\relax \ifnum 25<\count@\count@25\fi
  \def\x{\endgroup\@setsize\SetFigFont{#2pt}}%
  \expandafter\x
    \csname \romannumeral\the\count@ pt\expandafter\endcsname
    \csname @\romannumeral\the\count@ pt\endcsname
  \csname #3\endcsname}%
\fi
\fi\endgroup
\begin{picture}(10389,5922)(1784,-6142)
\put(2851,-1111){\makebox(0,0)[lb]{\smash{\SetFigFont{8}{9.6}{rm}limit group $D_\xi$}}}
\put(3751,-4786){\makebox(0,0)[lb]{\smash{\SetFigFont{8}{9.6}{rm}geodesic $L$}}}
\put(2926,-5836){\makebox(0,0)[lb]{\smash{\SetFigFont{8}{9.6}{rm}limit group $D_{-\xi}$}}}
\put(10351,-2161){\makebox(0,0)[lb]{\smash{\SetFigFont{8}{9.6}{rm}wandering sequence }}}
\put(9151,-6061){\makebox(0,0)[lb]{\smash{\SetFigFont{8}{9.6}{rm}limit group $K \ltimes U_-$}}}
\put(9151,-436){\makebox(0,0)[lb]{\smash{\SetFigFont{8}{9.6}{rm}limit group $K \ltimes U$}}}
\put(3826,-3736){\makebox(0,0)[lb]{\smash{\SetFigFont{8}{9.6}{rm}hyperbolic translation $\tau$}}}
\end{picture}
}
\pec

%
%

\title*{Appendix: Strong boundaries and commensurator super-rigidity}
\titlerunning{Commensurator super-rigidity}

\author{Patrick BONVIN}
\authorrunning{Patrick Bonvin}

\institute{}

\maketitle

\section*{Introduction}

Let $G$ be a locally compact second countable group, and let $\Gamma$
be a lattice of $G$, i.e.\ a
discrete subgroup such that $G/\Gamma$ carries a finite $G$-invariant measure.
The {\it commensurator~} of $\Gamma$ in $G$ is the group:
${\rm Comm}_G(\Gamma) := \{ g \! \in \! G \mid \Gamma \hbox{\rm ~and }
g\Gamma g^{-1} \hbox{\rm ~has finite
index in both }\Gamma \hbox{\rm ~and } g\Gamma g^{-1}\}$.

Our purpose is to use the recent double ergodicity theorem by V.
Kaimanovich on Poisson boundaries, in order to
show that G. Margulis' proof of the commensurator super-rigidity --
as analyzed and generalized by N. A'Campo and M.
Burger -- extends to a quite general setting.

\begin{thm} \label{main}
Let $G$ be a locally compact second countable topological group,
$\Gamma < G$ be a lattice and
$\Lambda$ be a subgroup of $G$ with $\Gamma < \Lambda < {\rm Comm}_G (\Gamma)$.
Let $k$ be a local field and $H$ be a connected almost $k$-simple
algebraic group.
Assume $\pi: \Lambda \to H_k$ is a homomorphism such that
$\pi(\Lambda)$ is Zariski dense in $H$ and
$\pi(\Gamma)$ is unbounded in the Hausdorff topology on $H_k$.
Then $\pi$ extends to a continuous homomorphism $\overline{\Lambda}
\to H_k / Z(H_k)$, where $Z(H_k)$ is the
center of $H_k$.
\end{thm}

This theorem is basically due to G. Margulis, who proved it in the
case where $G$ is a semisimple group over a locally
compact field \cite[VII.5.4]{Margulis}.
A deep idea in the proof is that the existence of the continuous
extension follows from the existence of a
$\Lambda$-equivariant map from the maximal Furstenberg boundary of
$G$ to a homogeneous space $H_k/L_k$,
where $L$ is a proper $k$-subgroup of $H$.
N. A'Campo and M. Burger extended the result to the case where $G$ is
as above \cite{AcaBur}, assuming the
existence of a closed subgroup $P$ playing the same measure-theoretic
role as a minimal parabolic subgroup.
This led M. Burger to state the above result assuming the existence
of a substitute for a maximal
Furstenberg boundary rather than a minimal parabolic subgroup \cite{Burger}.
The assumption that $k$ be of characteristic 0, made so far, was removed too.
M. Burger and N. Monod then constructed suitable boundaries for
compactly generated groups (up to
finite index, see \cite[Theorem 6]{BurMon}), which implied the above
result for a compactly generated group $G$
(\cite[Remark 7]{BurMon}).
The last step was made by V. Kaimanovich \cite{Kaim} (Theorem
\ref{spaceS} below), who proved that the Poisson
boundary for a nice measure on any  locally compact second countable
group is a strong hence a
suitable boundary.
We couldn't finish this historical summary without mentioning the
work of T.N. Venkataramana \cite{Venka}, who
was the first to prove super-rigidity and arithmeticity theorems in
arbitrary characteristics.

This note, which relies heavily on the proof given by N. A'Campo and
M. Burger in \cite{AcaBur},
shows how to use the previously quoted references to prove the above
commensurator super-rigidity.
It is organized as follows.
We first quote the results about boundaries of groups.
Then we recall how the existence of the continuous extension is
reduced to finding a $\Lambda$-equivariant
map from a boundary of $G$ to a non-trivial homogeneous space $H_k/L_k$.
We finally sketch the steps to construct the required
$\Lambda$-equivariant map, taking care of the fact that $k$ is
possibly of positive characteristic.

\section*{Poisson boundaries and strong boundaries}

Given a topological group $G$, a \textit{Banach G-module} is a pair
$(\pi,E)$ where $E$ is a Banach space and $\pi$
is an isometric linear representation of $G$ on $E$. The module
$(\pi,E)$ is \textit{continuous} if the action map
$G\times E \to E$ is continuous. A \textit{coefficient G-module} is a
Banach $G$-module $(\pi,E)$ contragredient to
some separable continuous Banach $G$-module, i.\ e.\ $E$ is the dual
of some separable Banach space $E^{\flat}$
and $\pi$ consists of operators adjoint to a continuous action of $G$
on $E^{\flat}$ (see \cite[chapter 1]{Monod}).
Denote by $\mathfrak{X}^{sep}$ the class of all separable coefficient modules.

Let $G$ be a locally compact group, and $(S,\mu)$ be a Lebesgue space
endowed with a measure class preserving
action of $G$. Given any class of coefficient Banach modules
$\mathfrak{X}$, the action of $G$ on $S$ is called
\textit{doubly $\mathfrak{X}$-ergodic} if for every coefficient
$G$-module $E$ in $\mathfrak{X}$, any weak-*
measurable $G$-equivariant function $f:S \times S \to E$ (with
respect to the diagonal action) is a.\ e.\ constant
(\cite[11.1.1]{Monod}).

Recall (\cite[4.3.1]{Zimmer}) that the $G$-action on $S$ is called
\textit{amenable} if for every separable Banach
space $E$ and every measurable right cocycle $\alpha: S \times G \to
{\rm Iso}(E)$ the following holds for
$\alpha^*$, the adjoint of the $\alpha$-twisted action on $L^1(S,E)$:
any $\alpha^*$-invariant measurable field
$\{A_s\}_{s \in  S}$ of non-empty convex weak-* compact subsets $A_s$
of the closed unit ball in $E^*$ admits a
measurable $\alpha^*$-invariant section.
We can now state V. Kaimanovich's result.

\begin{thm}\label{spaceS}
{\rm \cite{Kaim}} Let $G$ be a locally compact $\sigma$-compact group.
There exists a Lebesgue space $(S,\mu)$ endowed with a
measure class preserving action of $G$ such that:

\begin{itemize}
{\rm (i)~} The $G$-action on $S$ is amenable.

{\rm (ii)~}The $G$-action on $S$ is doubly $\mathfrak{X}^{sep}$-ergodic.
\end{itemize}
Such a space $(S, \mu)$ is called a \emph{strong $G$-boundary}
{\rm \cite[Def.\ 2.3]{MonSha}}.
\end{thm}

The space $S$ is a Poisson boundary for a suitable measure on $G$.
As mentioned before, this theorem strengthens a result of M. Burger
and N. Monod  \cite{BurMon}, who proved that any
compactly generated locally compact group possesses a finite index
open subgroup which has a strong boundary.

Note that the double $\mathfrak{X}^{sep}$-ergodicity of the
$G$-action on $S$ implies that the $\Gamma$-action
on $S$ is doubly $\mathfrak{X}^{sep}$-ergodic \cite[Prop.\
3.2.4]{BurMon} and that $G$ (as well as
any finite index subgroup of $\Gamma$)
acts ergodically on $S$ and on $S \times S$.

\section*{Reduction to finding a suitable equivariant map}

First, in view of the conclusion of the theorem, we may -- and shall
-- assume until the end of the note that the group
$H$ is adjoint.
Recall also that if $f:X\to Y$ is a measurable map from a Lebesgue
space $(X,\lambda)$ to a topological space $Y$,
its  \textit{essential image} is the closed subset of $Y$ defined by:
${\rm Essv}(f):= \{ y \! \in \! Y \mid \lambda( f^{-1}(V) ) > 0 \hbox{\rm ~for any
neighbourhood } V \hbox{\rm ~of } y \}$, and
$f$ is called \textit{essentially constant} if ${\rm Essv}(f)$
reduces to a point.
Here is the reduction theorem.

\begin{thm}\label{}
Let $G$ be a locally compact second countable topological group,
$\Gamma < G$ be a lattice and
$\Lambda$ be a subgroup of $G$ with $\Gamma < \Lambda < {\rm Comm}_G (\Gamma)$.
Let $k$ be a local field and $H$ be a connected almost $k$-simple
algebraic group.
Assume that $\pi: \Lambda \to H_k$ is a homomorphism such that
$\pi(\Lambda)$ is Zariski dense and that there
exists a $\Lambda$-equivariant non-essentially constant map $S \to
H_k/L_k$, where $L$ is a proper $k$-subgroup
of $H$.
Then there exists a continuous extension $\overline{\Lambda} \to H_k$ of $\pi$.
\end{thm}

The proof uses a simple and powerful ergodic argument \cite[Sect.
2.3]{AcaBur}, used many times in the full proof
of super-rigidity.
We will often deal with maps $\Theta : B \to M$ where $B$ is an
ergodic $\Gamma$-space, $M$ is a space with a
continuous $H_k$-action and $\Theta$ is equivariant with respect to a
group homomorphism $\Gamma \to H_k$.
Then if $M$ is a separable complete metrizable space and if the
$H_k$-orbits are locally closed in $M$, there is a
$H_k$-orbit $O$ in $M$ such that a conull subset of $B$ is sent to
$O$ by $\Theta$.
This is to be combined with the fact that a $k$-algebraic action of a
$k$-group $G$ on a $k$-variety $V$ induces a
continuous action of $G_k$ on $V_k$ in the Hausdorff topology, and
with the following crucial result, due to I.
Bernstein and A. Zelevinski.

\begin{thm}\label{closedorbits} {\rm \cite[6.15]{BerZel}}ð
Let k be a local field, V be a k-variety and G be a k-group acting
k-algebraically on V.
Then the orbits of $G_k$ in $V_k$ are locally closed.
\end{thm}

This theorem has a wide range of application because it implies local
closedness of orbits in many spaces.
Let $\mathcal{F}(S, W_k)$ be the space of classes of measurable maps
from $S$ to $W_k$,
endowed with the topology of convergence in measure on $\mathcal{F}(S, W_k)$.
It is metrizable by a complete separable metric.
Then, according to \cite[Lemma 6.7]{AcaBur}, the $H_k$-orbits in
$\mathcal{F}(S, W_k)$ are locally closed.
The proof there is given for a characteristic 0 local field $k$, but
it goes through once the
stabilizer ${\rm Stab}_H(w)$ of any $k$-rational point $w \! \in \!
W_k$, only $k$-closed in general, is replaced by
the $k$-subgroup $\overline{{\rm Stab}_{H_k}(w)}^Z$.
The ergodic argument applied to the function space $\mathcal{F}(S,
W_k)$ is a key point in the proof of the above
reduction theorem, whose proof can now be sketched (see \cite[Sect.
7]{AcaBur} for further details).

\textit{Proof.} Since $\overline{ \pi(\Lambda) }^Z = H$ and since $\theta$ is
$\Lambda$-equivariant,
we have $\overline{{\rm Essv}(\theta)}^Z=W$ where $W:=H/L$.
We define $\overline{\theta} : \overline{\Lambda} \to \mathcal{F}(S, W_k)$ by
$\overline{\theta}(\lambda)(s):=\theta(\lambda s)$.
It is $\Lambda$-equivariant and continuous.
Since $\Lambda$ acts ergodically on $\overline{\Lambda}$ and since
the $H_k$-orbits in $\mathcal{F}(S, W_k)$
are locally closed, there is a  $H_k$-orbit $O \subset
\mathcal{F}(S,W_k)$ such that $\overline{\theta}(\lambda)
\! \in \! O$ for almost all $\lambda \! \in \! \overline{\Lambda}$.
One deduces then from the fact that $O$ is open in $\overline{O}$ and
$\overline{\theta}$ is continuous, that
$\overline{\theta}(\overline{\Lambda}) \subset O$.
In particular $O = {H_k}_* \theta$.
Then it follows from $\overline{{\rm Essv}(\theta)}^Z = W$ that
${\rm Stab}_{H_k} (\theta)$ fixes pointwise $W$
and thus is trivial since $H$ is adjoint.
Therefore the map $h : \overline{ \Lambda } \to H_k$ defined by
$\overline{\theta} (\lambda) = h(\lambda)_* \theta$ for any $\lambda
\! \in \! \overline{\Lambda}$, is  a
continuous homomorphism.
Since $\theta$ is $\Lambda$-equivariant, $h$ is the desired extension of $\pi$. \carre

\section*{Constructing the required equivariant map}

We now sketch the proof of the existence of a $\Lambda$-equivariant
map as above under the hypotheses of
Theorem \ref{main}.
Since $k$ is of arbitrary characteristic, the adjoint representation
Ad of $H$ need no longer be irreducible.
Still, we can choose $\rho : H \to {\rm GL}(V)$ a faithful rational
representation of $H$, defined and irreducible
over $k$, on a finite-dimensional $k$-vector space $V$.
The induced map $\rho : H_k \to {\rm PGL}(V_k)$ is injective because
$H$ is adjoint, and by \cite[I.2.1.3]{Margulis}
it is a closed embedding.
We have a homomorphism $\rho \circ \pi : \Gamma \to {\rm PGL}(V_k)$,
so that  $\Gamma$ acts by
homeomorphisms on the compact metric space $\prv_k$.
This induces a continuous action $\Gamma \times M^1(\prv_k) \to
M^1(\prv_k)$, where $M^1(\prv_k)$ is the
space of probability measures on $\prv_k$ endowed with the compact
metrizable weak-$*$ topology.

\begin{prop}\label{equiv}
Let $G$ be a locally compact group and $(S,\mu)$ be a Lebesgue space
on which $G$ acts amenably.
Then, possibly after discarding an invariant null set in $S$, there 
exists a measurable
$\Gamma$-equivariant map
$\phi: S \to M^1(\prv_k)$.
\end{prop}

\textit{Proof. } This follows immediately from Theorem 4.3.5 and Proposition 4.3.9 in
\cite{Zimmer}. \carre

At this stage, we have a measurable map $\phi: S \to M^1(\prv_k)$
which is equivariant for the $\Gamma$-action
only, and which goes to a space of probability measures.
The next step provides a $\Gamma$-equivariant map to a homogeneous
space $H_k/L_k$.

We denote by
${\rm Var}_k (\prv)$ the set of algebraic subvarieties of $\prv$
defined over $k$ and by
${\rm supp}_Z : M^1({\bf P} V_k) \to {\rm Var}_k (\prv)$  the map
which to a probability measure $\mu$
associates $\overline{{\rm supp}(\mu)}^Z$, the Zariski closure of its support.
For any $n$-dimensional $k$-vector space $W_k$ we set
${\rm Gr}(W_k)  := \bigsqcup_{l=0}^n {\rm Gr}_l(W_k)$,
where ${\rm Gr}_l(W_k)$ is the compact Grassmannian of $l$-planes in $W_k$.
By attaching to each projective variety $X \subset \prv$ its graded
defining ideal $I_X$, we see
${\rm Var}_k (\prv)$ as a subspace of the compact space
$\prod_{d=0}^\infty {\rm Gr}(k[V]_d)$, where $k[V]_d$ is the space of
$d$-homogeneous polynomials on $V_k$.
This induces a topology on ${\rm Var}_k (\prv)$, and it is proved in
\cite[Sect. 5]{AcaBur}, by
characteristic free arguments, that the map ${\rm supp}_Z$ is
measurable and ${\rm PGL}(V_k)$-equivariant.
Therefore we obtain by composition a $\Gamma$-equivariant measurable
map $\Phi: S \to {\rm Var}_k (\prv)$,
sending each $s \! \in \! S$ to the Zariski closure of the support
of $\phi(s)$.
We denote it by $\Phi$, and call it \textit{boundary map}.

\begin{thm}\label{nonconstant}
{\rm \cite[Theorem 5.1]{AcaBur}}ð
The boundary map $\Phi$ is not essentially constant.
\end{thm}

This result follows from the arguments in \cite[Sect. 5]{AcaBur}.
To see this, we first note that since $H$ is $k$-simple,
$\pi(\Gamma)$ is unbounded and $\pi(\Lambda)$ is
Zariski dense, the inclusion $\Lambda < {\rm Comm}_G(\Gamma)$ and
the fact that the identity
component of an algebraic group is always a finite index subgroup
imply that $\pi(\Gamma)$ is Zariski dense in $H$.
The other facts needed in \cite[Sect. 5]{AcaBur} are the ergodicity
of $\Gamma$  on $S$ and on
$S \times S$, and the Furstenberg lemma, all  available in our context.

  From the result, there is a $d$ for which $\Phi : S \to {\rm
Gr}(k[V]_d)$ is not essentially constant.
The ergodic argument of the previous section and the ergodicity  of
the $\Gamma$-action on $S$ imply that $\Phi$
essentially sends $S$ to a $H_k$-orbit in ${\rm Gr}(k[V]_d)$, which
is homeomorphic to a space $H_k / L_k$ for
some proper algebraic subgroup $L$ of $H$: we have obtained a
$\Gamma$-equivariant measurable map
$\phi: S \to H_k / L_k$.

The very last step consists in passing from $\Gamma$- to
$\Lambda$-equivariance.
Once maps as above are known to exist, the descending chain condition
for algebraic sugbroups and Zorn's lemma, as
used in \cite[Sect.7]{AcaBur}, prove the existence of a couple
$(\phi, H_k / L_k)$ satisfying a universal property.
The normalizer of $L_k$ in $H$ may only be $k$-closed, but if we
denote by $L'$ the Zariski closure of the
normalizer of $L_k$ in $H_k$, we get a $k$-subgroup, which is proper
by $k$-simplicity of $H$ and such that:

\begin{thm}\label{thetaequ}
{\rm \cite[Corollary 7.2]{AcaBur}}ð
The composed map $\theta : S \to H_k / L_k \to H_k / L_k'$ is 
$\Lambda$-equivariant
and measurable.
\end{thm}

\bibliographystyle{amsalpha}
\bibliography{csr}

\vspace{.5cm}

\textsc{ETHZ -- Department of Mathematics\\
R\"amistrasse 101 \\ 
8092 Z\"urich -- Switzerland} \\
\textit{E-mail address:} \texttt{bonvin@math.ethz.ch}

\end{document}